\documentclass{article}
\usepackage[utf8]{inputenc}
\usepackage[T1]{fontenc}
\usepackage{lmodern}
\usepackage{babel}
\usepackage{lipsum}
\usepackage{graphicx}
\usepackage{amsmath}
\usepackage{amsfonts}
\usepackage{subcaption}
\usepackage{tabularx}
\usepackage{siunitx}
\usepackage{booktabs}
\usepackage{dashrule}
\usepackage{tikz}
\usetikzlibrary{calc}
\usetikzlibrary{patterns}
\usepackage[ruled,linesnumbered]{algorithm2e}
\usepackage{amsthm}
\usepackage{accents}
\usepackage{parskip}
\usepackage{amssymb}
\usepackage[hyperfootnotes=false]{hyperref}
\usepackage{thmtools}
\usepackage{multirow}
\usepackage{alphabeta} % Allows \alpha, etc. in text mode.
\usepackage{thm-restate}

\newtheorem{definition}{Definition}

\usepackage{todonotes}
\usepackage[scaled=0.90]{FiraMono}

\usepackage{arxiv}
\usepackage{charter}
\usepackage{fullpage}

\hypersetup{
    colorlinks=true,
    linkcolor=blue,
    filecolor=magenta,      
    urlcolor=cyan,
    citecolor=blue
}

\newcommand{\coefficients}{\ensuremath{\boldsymbol{\alpha}}\xspace}

\newcommand{\lprelaxation}{\ensuremath{\mathcal{P}}\xspace}

\newcommand{\xv}{\ensuremath{\mathbf{x}}\xspace}
\newcommand{\lpoptimal}{\ensuremath{\mathbf{x}^{LP}}\xspace}
\newcommand{\lpoptimali}[1]{\ensuremath{x^{LP}_{#1}}\xspace}
\newcommand{\optimal}{\ensuremath{\mathbf{x}^{*}}\xspace}

\newcommand{\incumbent}{\ensuremath{\hat{\mathbf{x}}}\xspace}
\newcommand{\basissol}{\ensuremath{\bar{\mathbf{x}}}\xspace}
\newcommand{\basissoli}[1]{\ensuremath{\bar{x}_{#1}}\xspace}

\newcommand{\cut}{\ensuremath{(\boldsymbol{\alpha},\beta)}\xspace}

\newcommand{\reals}{\ensuremath{\mathbb{R}}\xspace}

\newcommand{\ints}{\ensuremath{\mathbb{Z}}\xspace}

\newcommand{\sset}{\ensuremath{\mathcal{S}}\xspace}
\newcommand{\disjunction}{\ensuremath{\mathcal{D}}\xspace}
\newcommand{\disjunctioni}[1]{\ensuremath{\mathcal{D}_{#1}}\xspace}
\newcommand{\splitcoeff}{\ensuremath{\boldsymbol{\pi}}\xspace}
\newcommand{\splitrhs}{\ensuremath{\pi_{0}}\xspace}
\newcommand{\splitdisjunction}[2]{\ensuremath{\mathcal{D}(\boldsymbol{\pi},\pi_{0})}\xspace}
\newcommand{\lpcone}[1]{\ensuremath{\mathcal{C}(#1)}\xspace}

\newcommand{\eff}[3]{\ensuremath{\mathtt{eff}(#1,#2,#3)}\xspace}

\newcommand{\floor}[1]{\lfloor #1 \rfloor}
\newcommand{\ciel}[1]{\lceil #1 \rceil}

\title{Branching via Cutting Plane Selection:\\ Improving Hybrid Branching}

\author{ \href{https://orcid.org/0000-0001-7270-1496}{\includegraphics[scale=0.06]{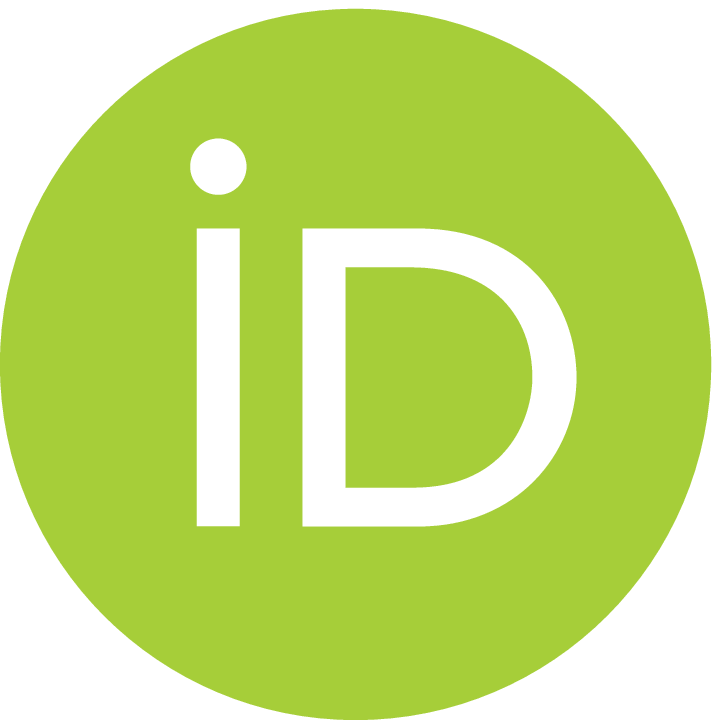}}\hspace{1mm}Mark Turner\thanks{Institute of Mathematics, Technische Universit{\"a}t Berlin, Straße des 17. Juni 135, 10623 Berlin, Germany}\hspace{2mm}\thanks{Zuse Institute Berlin, Department of Mathematical Optimization, Takustr. 7, 14195 Berlin} \\
	\texttt{turner@zib.de} \\
	\And
	\href{https://orcid.org/0000-0002-6320-8154}{\includegraphics[scale=0.06]{orcid_id_icon.eps}}\hspace{1mm}Timo Berthold\footnotemark[1]\hspace{2mm}\thanks{Fair Isaac Germany GmbH, Takustr. 7, 14195 Berlin, Germany} \\
	\texttt{timoberthold@fico.com} \\
	\And
	\href{https://orcid.org/0000-0002-6284-3033}{\includegraphics[scale=0.06]{orcid_id_icon.eps}}\hspace{1mm}Mathieu Besançon\footnotemark[2]\\
	\texttt{besancon@zib.de} \\
	\And
	\href{https://orcid.org/0000-0002-1967-0077}{\includegraphics[scale=0.06]{orcid_id_icon.eps}}\hspace{1mm}Thorsten Koch\footnotemark[1]\hspace{2mm}\footnotemark[2]
	\\
	\texttt{koch@zib.de}
}

\begin{document}

% \ZTPAuthor{\ZTPHasOrcid{Mark Turner}{0000-0001-7270-1496},
% \ZTPHasOrcid{Timo Berthold }{0000-0002-6320-8154},
% \ZTPHasOrcid{Mathieu Besançon}{0000-0002-6284-3033}, \\
% \ZTPHasOrcid{Thorsten Koch }{0000-0002-1967-0077}}
% \ZTPTitle{Branching via Cutting Plane Selection:\\ Improving Hybrid Branching}
% \ZTPNumber{23-17}
% \ZTPMonth{June}
% \ZTPYear{2023}

% \zibtitlepage

\maketitle

\begin{abstract}
Cutting planes and branching are two of the most important algorithms for solving mixed-integer linear programs. For both algorithms, disjunctions play an important role, being used both as branching candidates and as the foundation for some cutting planes. We relate branching decisions and cutting planes to each other through the underlying disjunctions that they are based on, with a focus on Gomory mixed-integer cuts and their corresponding split disjunctions. We show that selecting branching decisions based on quality measures of Gomory mixed-integer cuts leads to relatively small branch-and-bound trees, and that the result improves when using cuts that more accurately represent the branching decisions. Finally, we show how the history of previously computed Gomory mixed-integer cuts can be used to improve the performance of the state-of-the-art hybrid branching rule of SCIP. Our results show a $4\%$ decrease in solve time, and an $8\%$ decrease in number of nodes over affected instances of MIPLIB 2017. 
\end{abstract}

\section{Introduction}

This paper proposes a new criterion for branching in branch-and-cut algorithms to solved Mixed-Integer Linear Programs (MILP) based on the disjunctions that underpin both branching candidates and several families of cutting planes. A MILP is an optimisation problem that is classically defined as:
\begin{align}
    \underset{\mathbf{x}}{\mathrm{min}}\{\mathbf{c}^{\intercal}\mathbf{x} \;\; | \;\; \mathbf{A}\mathbf{x} \leq \mathbf{b}, \;\; \mathbf{l} \leq \mathbf{x} \leq \mathbf{u}, \;\; \mathbf{x} \in \mathbb{Z}^{|\mathcal{J}|} \times \mathbb{R}^{n - |\mathcal{J}|} \} \label{eq:mip}
\end{align}
Here, $\mathbf{c} \in \reals^{n}$ is the objective coefficient vector, $\mathbf{A} \in \reals^{m \times n}$ is the constraint matrix, $\mathbf{b} \in \reals^{m}$ is the right hand side constraint vector, $\mathbf{l}, \mathbf{u} \in \reals^{n} \cup \{-\infty, \infty\}^{n}$ are the lower and upper variable bound vectors, and $\mathcal{J} \subseteq \{1 , \dots , n\}$ is the set of indices of integer variables.
We denote the the feasible region of the linear programming (LP) relaxation of \eqref{eq:mip} as \lprelaxation, where the LP is derived by relaxing the integrality requirements of \eqref{eq:mip}. 
An optimal solution to \eqref{eq:mip} is denoted \optimal, a feasible solution denoted \incumbent, and an LP optimal solution is denoted as \lpoptimal.

\begin{figure}[h]
\centering
\begin{subfigure}[b]{0.475\textwidth}
\centering
    \begin{tikzpicture}[scale=1]
    % Create the grid points
    % Create the grid points
    \foreach \x/\y in {0/0, 0/1, 0/2, 0/3, 1/0, 1/1, 1/2, 1/3, 2/0, 2/1, 2/2, 2/3, 3/0, 3/1, 3/2, 3/3}
    {
    \fill[opacity=0.4] (\x,\y) circle (2pt);
    }
    \foreach \x/\y in {1/0, 2/0, 2/1, 3/1, 3/2}
    {
    \fill[opacity=1] (\x,\y) circle (2pt);
    }
    % Create the nodes of the polygon
    \node (a) at (1.5,1.5) {};
    \node (b) at (3,2.5) {};
    \node (c) at (3,1) {};
    \node (d) at (2.5,0) {};
    \node (e) at (0.5,0) {};
    
    \node (d11) at (1,0.75) {};
    \node (d12) at (1,0) {};
    \node (d21) at (2,1.8333) {};
    \node (d22) at (2,0) {};
    
    % draw the big polygon    
    \draw[thick] (a.center) -- (b.center) -- (c.center) -- (d.center)  -- (e.center) -- cycle;
    % Draw the smaller disjunctions
    \filldraw[thick,fill=blue!60,fill opacity=0.4] (d11.center) -- (d12.center) -- (e.center) -- cycle;
    \filldraw[thick,fill=blue!60,fill opacity=0.4] (d21.center) -- (b.center) -- (c.center) -- (d.center) -- (d22.center) -- cycle;
    % Label the LP optimal point
    \fill[draw=black,fill=red] (a.center) circle (2pt);
    
    % Draw the S free set
    \draw[thick,dashed] (1, -0.5) -- (1, 3.5);
    \draw[thick,dashed] (2, -0.5) -- (2, 3.5);
    
    \end{tikzpicture}
\label{fig:branching_1}
\end{subfigure}
\hfill
\begin{subfigure}[b]{0.475\textwidth}
\centering
    \begin{tikzpicture}[scale=1]
    % Create the grid points
    \foreach \x/\y in {0/0, 0/1, 0/2, 0/3, 1/0, 1/1, 1/2, 1/3, 2/0, 2/1, 2/2, 2/3, 3/0, 3/1, 3/2, 3/3}
    {
    \fill[opacity=0.4] (\x,\y) circle (2pt);
    }
    \foreach \x/\y in {1/0, 2/0, 2/1, 3/1, 3/2}
    {
    \fill[opacity=1] (\x,\y) circle (2pt);
    }
    % Create the nodes of the polygon
    \node (a) at (1.5,1.5) {};
    \node (b) at (3,2.5) {};
    \node (c) at (3,1) {};
    \node (d) at (2.5,0) {};
    \node (e) at (0.5,0) {};
    
    % draw the big polygon    
    \draw[thick] (a.center) -- (b.center) -- (c.center) -- (d.center)  -- (e.center) -- cycle;
    % Draw the smaller polygon
    \draw[thick,fill=blue!60, fill opacity=0.4] (2.25, 2) -- (b.center) -- (c.center) -- (d.center) -- (e.center) -- (1, 0.75) -- cycle;

    % Label the LP optimal point
    \fill[draw=black,fill=red] (a.center) circle (2pt);

    \draw[ultra thick, red] (-0.25, -0.5) -- (3.5, 3.25);
    
    \end{tikzpicture}
\label{fig:branching_2}
\end{subfigure}
\caption{(Left) An example branching decision. The red point is the LP optimal solution, the larger polytope the original LP feasible region, and the two blue polytopes the feasible LP regions of the subproblems. (Right) An example cutting plane.}
\label{fig:branching}
\end{figure}
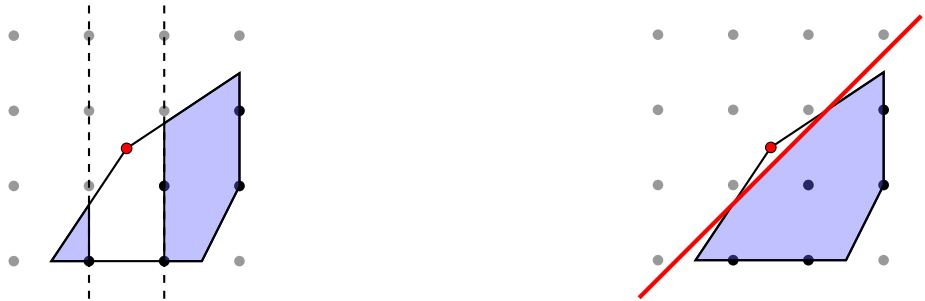

The core algorithm for solving MILPs is \emph{branch-and-cut}, see \cite{achterberg2007constraint} for a thorough introduction of MILP solving methods.
Branch-and-cut combines two main components, \emph{branch-and-bound} and \emph{cutting planes}.
The branch-and-bound algorithm recursively partitions the MILP search space
by a process called \emph{branching} i.e.~splitting a problem into smaller subproblems. This recursion creates a \emph{tree}, where each node is a subproblem.
Traditionally branching is performed on an integer variable, $x_{i}$ with fractional LP value, $x^{LP}_{i}$, creating the two LP subproblems with feasible regions $\lprelaxation\, \cap\, \{x^{LP}_{i} \leq \floor{x^{LP}_{i}}\}$ and $\lprelaxation\, \cap \, \{x^{LP}_{i} \leq \ciel{x^{LP}_{i}}\}$, thus making the current LP solution infeasible in both children problems.
An example branching procedure is visualised in Figure~\ref{fig:branching}.
The algorithm bounds the optimal objective through upper bounds from feasible solutions of Problem~\eqref{eq:mip} obtained at leaf nodes of the tree, and lower bounds from LP relaxations that bound all subtrees of a node.
The algorithm of branch-and-bound that we focus on is \emph{variable selection}, which is concerned with determining which variable to branch on at a given node from the given candidates. 

A cutting plane, or \emph{cut}, parameterised by $(\coefficients, \beta) \in \reals^{n+1}$, is an inequality $\coefficients^{\intercal} \xv \leq \beta$ that is violated by at least one solution of the LP relaxation but that does not increase the optimal value of the problem when added to the problem, i.e., it is valid for \eqref{eq:mip}.
This definition of a cut is more general than the classical one, which requires that a cut does not remove any integer-feasible solution to \eqref{eq:mip}, and separates some LP feasible fractional solution.
Our definition however captures additional families of cuts such as symmetry-breaking cuts \cite{pfetsch2019computational,hojny2019polytopes}.
The cutting plane algorithm iteratively generates cuts, applies them, and re-solves the tightened LP relaxation. In the classical case, this algorithm is repeated until an integral LP relaxation solution is achieved. In branch-and-cut, this algorithm is repeated at the root node until some termination criteria is met, with additional cuts applied throughout the branch-and-bound tree. The procedure of cutting planes that we focus on is \emph{cut selection}, which is concerned with deciding which subset of the computed cuts to actually add to the LP relaxation. 

In this paper we propose a technique for the variable selection problem based on cut selection. We note that the opposite direction -- using variable selection techniques for cut selection -- also presents a valid avenue of research, however this lies outside the scope of this paper.
For each branching candidate, we generate a cut, greedily select the best cut according to standard cut scoring measures, and then branch on the corresponding candidate. Specifically, we will generate Gomory Mixed-Integer (GMI) cuts, for which we provide a detailed introduction in Section~\ref{sec:gmi}. This approach to variable selection is objective-free, i.e.~not based on the objective vector, and is thus complementary to standard pseudo-cost based approaches \cite{achterberg2005branching}.
We show the effectiveness of multiple variants of this branching rule using different levels of strengthened cuts, and compare them to standard branching rules from the literature. Finally, we show how this information can be incorporated into the hybrid branching rule of SCIP~\cite{achterberg2009hybrid,scip8}, resulting in an improvement of general solver performance.

\section{Related Work}

Branching in MILP has been thoroughly studied, both theoretically and computationally, see \cite{linderoth1999computational,achterberg2005branching,dey2021theoretical}. The current state-of-the-art variable selection method is hybrid branching \cite{achterberg2009hybrid}, which is reliability pseudo-cost branching \cite{achterberg2005branching} with integrated constraint satisfaction and satisfiability problem techniques. An array of other selection rules exist, such as nonchimerical fractionality branching \cite{fischetti2012branching}, cloud branching \cite{berthold2013cloud}, and general disjunction branching \cite{owen2001experimental}. The above-stated methods, unlike our cut selection approach to branching, often depend on the objective function for variable selection.
We note that variable selection has also served as the playground for introducing machine learning to MILP solvers, see \cite{khalil2016learning,gasse2019exact} for early examples, and \cite{huang2021branch} for an overview.

Cutting plane selection has been less studied than variable selection. It is however currently experiencing a recent refocus through machine learning-driven research, see \cite{deza2023machine} for an overview.
Early computational studies, see \cite{achterberg2007constraint,wesselmann2012implementing}, show that a diverse set of measures is necessary for good performance when scoring cuts. Both studies, as well as the computational study \cite{zerohalf}, use parallelism-based cut filtering algorithms, and show that the inclusion of filtering methods is critical to performance. These studies suggest that it is preferable for performance to select from a large set of weaker cuts than from a small set of stronger cuts. More recent work on cut selection, for which \cite{dey2018theoretical} provides ample motivation, is machine learning-based.
Specifically, research has focused on theoretical guarantees \cite{balcan2021sample,balcan2022improved,turner2022adaptive}, new scoring measures \cite{turner2022cutting}, the amount of cuts to select \cite{wang2023learning}, and learning to score cuts with supervised \cite{baltean2019scoring,huang2021learning}, imitation \cite{paulus2022learning}, and reinforcement learning \cite{tang2020reinforcement,wang2023learning}.

Our work is not the first to use cutting plane selection to dictate branching decisions.
Moreover, it is not the first to use GMI cuts specifically, see \cite{karamanov2011branching,cornuejols2011improved}.
In both papers, the split disjunctions, which define the GMI cuts of tableau rows, are used as branching candidates. The efficacy of the GMI cuts are used to filter the set of branching candidates, where ultimately strong branching is used as the final selection criterion. In \cite{cornuejols2011improved}, additional experiments are presented that compare disjunctions derived from reduce-and-split cuts, see \cite{andersen2005reduce}.
Our research differs from \cite{karamanov2011branching,cornuejols2011improved} in that we branch on elementary splits, i.e., single variable disjunctions, we perform extensive computational experiments using non-strengthened GMI cuts, and we integrate our approach within existing state-of-the-art history-based methods in a modern MILP solver.

\section{Gomory Mixed-Integer Cuts}
\label{sec:gmi}

This section provides a thorough introduction to Gomory Mixed-Integer (GMI) cuts.
Following the history and general introduction of Subsection~\ref{subsec:gmi_history}, we introduce disjunctive, split, and intersection cuts in Subsection~\ref{subsec:intersection_and_split}. We then step through the derivation of GMI inequalities in Subsection~\ref{subsec:gmi_derive}, ending with how GMI cuts are used and derived in practice in Subsection~\ref{subsec:gmi_practice}. The geometric interpretation of GMI cuts will be provided, and linked to the relations between the families of cuts introduced in Subsection \ref{subsec:intersection_and_split}.
For alternative overviews of GMI cuts, see \cite{jorg2008k,cornuejols2008valid,andersen2005reduce,cornuejols2011improved}.

\subsection{GMI Introduction and History}
\label{subsec:gmi_history}

First introduced in 1960~\cite{gomory1960algorithm}, GMI inequalities are general-purpose inequalities valid for arbitrary bounded MILPs. They can be used to iteratively tighten an LP relaxation of an MILP, and when they are generated to separate a specific solution, are referred to as cutting planes, or \textit{cuts}. In practice, they are generated to separate the current LP solution using the simplex tableau, see Subsection~\ref{subsec:gmi_practice}. 

Following the landmark paper \cite{balas1996gomory}, GMI cuts were empirically shown to be a computational success. This success was in spite of a commonly held belief that only cuts derived from the specific structure of an MILP instance were computationally useful. Some examples of structured inequalities or cuts are knapsack cover and flow cover inequalities \cite{conforti2014integer}. The summarised reasons for the success of \cite{balas1996gomory} was their intelligent lifting procedure to globally valid cuts, their selection algorithm, their use of branch-and-cut as opposed to pure cutting plane approaches, and the increased robustness of LP solvers. For a more complete history behind the resurgence of GMI inequalities, see \cite{cornuejols2007revival}. Advances on GMI cuts have continued, where we name reduce-and-split cuts \cite{andersen2005reduce} and LaGromory cuts \cite{fischetti2011relax} as examples. To stress the importance of these cuts in the current MILP landscape, we note that GMI cuts are continually noted as computationally necessary \cite{wesselmann2009strengthening,achterberg2013mixed}, and are used in every state-of-the-art MILP solver, see Xpress \cite{xpress}, Gurobi \cite{gurobi}, CPLEX \cite{cplex}, HiGHS \cite{highs}, and SCIP \cite{scip8}.

\subsection{Disjunctive, Split, and Intersection Cuts}
\label{subsec:intersection_and_split}

It is common in the literature to find compact introductions of GMI cuts that mention they are either disjunctive cuts, intersection cuts, or split cuts. All these statements are true, and moreover, the families of cuts have a clear hierarchy \cite{balas2018disjunctive, andersen2005split}:
\begin{align*}
    \text{GMI cuts from basic feasible solutions} \; \subset \; \text{Split cuts} \; \subset \; \text{Intersection cuts} \; \subset \; \text{Disjunctive cuts}
\end{align*}

We highlight that while our work on branching leverages GMI cuts, it can also leverage any family of cuts that fit into this hierarchy and can be derived from disjunctions.

\subsubsection{Disjunctions and Disjunctive Cuts}

A linear disjunction is a set of linear inequalities joined by \emph{and}, \emph{or}, and \emph{negation} operators (see \cite{balas2018disjunctive} for a thorough introduction). The solution set of a disjunction is a \emph{disjunctive set}. For MILPs, every integer-feasible solution to \eqref{eq:mip} is an element of a disjunctive set. A \emph{disjunctive cut} is any cut derived from such a disjunctive set, i.e., all elements in the disjunctive set remain feasible and some fractional solution outside the disjunctive set is separated.

\begin{figure}[h]
\centering
\begin{subfigure}[b]{0.475\textwidth}
    \centering
    \begin{tikzpicture}[scale=1]
    % Create the grid points
    % Create the grid points
    \foreach \x/\y in {0/0, 0/1, 0/2, 0/3, 1/0, 1/1, 1/2, 1/3, 2/0, 2/1, 2/2, 2/3, 3/0, 3/1, 3/2, 3/3}
    {
    \fill[opacity=0.4] (\x,\y) circle (2pt);
    }
    \foreach \x/\y in {1/0, 2/0, 2/1, 3/1, 3/2}
    {
    \fill[opacity=1] (\x,\y) circle (2pt);
    }
    % Create the nodes of the polygon
    \node (a) at (1.5,1.5) {};
    \node (b) at (3,2.5) {};
    \node (c) at (3,1) {};
    \node (d) at (2.5,0) {};
    \node (e) at (0.5,0) {};
    
    \node (d21) at (1,0.75) {};
    \node (d22) at (1,0) {};
    \node (d31) at (2.25,2) {};
    \node (d32) at (3,2) {};
    \node (d41) at (2,0) {};
    \node (d42) at (2,1) {};
    
    % draw the big polygon    
    \draw[thick] (a.center) -- (b.center) -- (c.center) -- (d.center)  -- (e.center) -- cycle;
    % Draw the smaller disjunctions
    \filldraw[thick,fill=blue!60,fill opacity=0.4] (d21.center) -- (d22.center) -- (e.center) -- cycle;
    \filldraw[thick,fill=blue!60,fill opacity=0.4] (d31.center) -- (b.center) -- (d32.center) -- cycle;
    \filldraw[thick,fill=blue!60,fill opacity=0.4] (c.center) -- (d.center) -- (d41.center) -- (d42.center) -- cycle;
    % Label the LP optimal point
    \fill[draw=black,fill=red] (a.center) circle (2pt);
    
    % Draw the S free set
    \draw[thick,dashed] (-0.5,1) -- (1,1);
    \draw[thick,dashed] (2,1) -- (3.5,1);
    \draw[thick,dashed] (-0.5,2) -- (1,2);
    \draw[thick,dashed] (2,2) -- (3.5,2);
    \draw[thick,dashed] (1,3.5) -- (1,2);
    \draw[thick,dashed] (1,-0.5) -- (1,1);
    \draw[thick,dashed] (2,3.5) -- (2,2);
    \draw[thick,dashed] (2,-0.5) -- (2,1);
    \node (f) at (1.5,1.8) {};
    
    \node[label=above:{\disjunctioni{1}}] (d1) at (0.5,2) {};
    \node[label=above:{\disjunctioni{2}}] (d2) at (0.5,0.25) {};
    \node[label=above:{\disjunctioni{3}}] (d3) at (2.5,2.25) {};
    \node[label=above:{\disjunctioni{4}}] (d4) at (2.5,0.25) {};
    
    \end{tikzpicture}
\label{fig:disjunction_1}
\end{subfigure}
\hfill
\begin{subfigure}[b]{0.475\textwidth}
\centering
    \begin{tikzpicture}[scale=1]
    % Create the grid points
    \foreach \x/\y in {0/0, 0/1, 0/2, 0/3, 1/0, 1/1, 1/2, 1/3, 2/0, 2/1, 2/2, 2/3, 3/0, 3/1, 3/2, 3/3}
    {
    \fill[opacity=0.4] (\x,\y) circle (2pt);
    }
    \foreach \x/\y in {1/0, 2/0, 2/1, 3/1, 3/2}
    {
    \fill[opacity=1] (\x,\y) circle (2pt);
    }
    % Create the nodes of the polygon
    \node (a) at (1.5,1.5) {};
    \node (b) at (3,2.5) {};
    \node (c) at (3,1) {};
    \node (d) at (2.5,0) {};
    \node (e) at (0.5,0) {};
    
    \node (d21) at (1,0.75) {};
    \node (d22) at (1,0) {};
    \node (d31) at (2.25,2) {};
    \node (d32) at (3,2) {};
    \node (d41) at (2,0) {};
    \node (d42) at (2,1) {};
    
    % draw the big polygon    
    \draw[thick] (a.center) -- (b.center) -- (c.center) -- (d.center)  -- (e.center) -- cycle;
    % Draw the smaller disjunctions
    \filldraw[thick,fill=blue!60,fill opacity=0.4] (d21.center) -- (d22.center) -- (e.center) -- cycle;
    \filldraw[thick,fill=blue!60,fill opacity=0.4] (d31.center) -- (b.center) -- (d32.center) -- cycle;
    \filldraw[thick,fill=blue!60,fill opacity=0.4] (c.center) -- (d.center) -- (d41.center) -- (d42.center) -- cycle;
    % Label the LP optimal point
    \fill[draw=black,fill=red] (a.center) circle (2pt);
    
    % Draw the S free set
    \draw[thick,dashed] (-0.5,1) -- (1,1);
    \draw[thick,dashed] (2,1) -- (3.5,1);
    \draw[thick,dashed] (-0.5,2) -- (1,2);
    \draw[thick,dashed] (2,2) -- (3.5,2);
    \draw[thick,dashed] (1,3.5) -- (1,2);
    \draw[thick,dashed] (1,-0.5) -- (1,1);
    \draw[thick,dashed] (2,3.5) -- (2,2);
    \draw[thick,dashed] (2,-0.5) -- (2,1);
    \node[] (f) at (1.5,1.8) {};
    
    % \node[label=above:{\disjunctioni{1}}] (d1) at (0.5,2) {};
    % \node[label=above:{\disjunctioni{2}}] (d2) at (0.5,0.25) {};
    % \node[label=above:{\disjunctioni{3}}] (d3) at (2.5,2.25) {};
    % \node[label=above:{\disjunctioni{4}}] (d4) at (2.5,0.25) {};

    \draw[ultra thick, red] (-0.25, -0.5) -- (3.5, 3.25);
    
    \end{tikzpicture}
\label{fig:disjunction_2}
\end{subfigure}
\caption{(Left) An example disjunction $((x_{1} \leq \floor{\lpoptimali{1}}) \lor (x_{1} \geq \ciel{\lpoptimali{1}})) \land ((x_{2} \leq \floor{\lpoptimali{2}}) \lor (x_{2} \geq \ciel{\lpoptimali{2}}))$ . The disjunctive set is the union of blue polytopes. Here \disjunctioni{1} is the empty set. (Right) An example disjunctive cut for the disjunction.}
\label{fig:disjunction}
\end{figure}
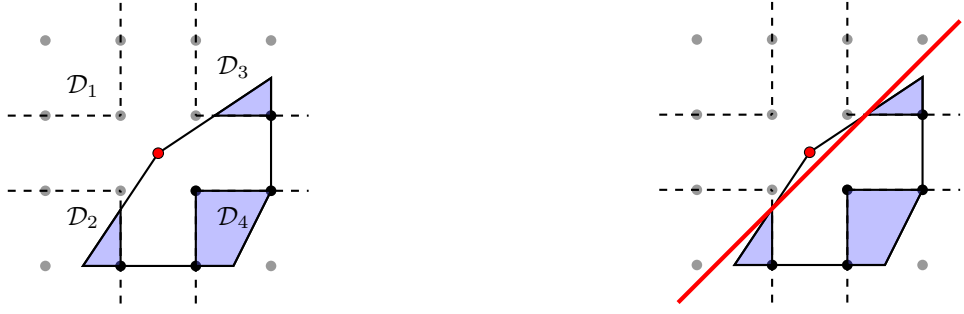

A linear disjunctive set represents a union of polyhedra. It is defined as:
% \bigcup_{i \in \{1,...,|\disjunction|\}}
\begin{align*}
    \disjunction := \bigcup_{i=1}^{|\disjunction|} \disjunctioni{i}, \;\; \text{where} ; \disjunctioni{i} \subseteq \lprelaxation \;\; \forall i \in \{1, \cdots |\disjunction|\}
\end{align*}
An example disjunction is visualised in Figure~\ref{fig:disjunction}, along with a valid disjunctive cut.

\subsubsection{Splits and Split Cuts}

A \emph{split disjunction}, or \emph{split}, is defined by an integer $\splitrhs \in \ints$ and an integral vector $\splitcoeff \in \ints^{|\mathcal{J}|} \times \mathbf{0}^{n - |\mathcal{J}|}$, which has zero entries for coefficients of continuous variables. We denote the split disjunction as \splitdisjunction{\splitcoeff}{\splitrhs}, where $(\splitcoeff, \splitrhs)$ define the two hyperplanes:
\begin{align}
\begin{split}
    \splitcoeff^{\intercal} \mathbf{x} &\leq \splitrhs \\ \splitcoeff^{\intercal} \mathbf{x} &\geq \splitrhs + 1 \label{eq:split}
\end{split}
\end{align}
The disjunctive set $\disjunction = \bigcup_{i \in \{1,2\}} \disjunctioni{i}$ formed by the hyperplanes is:
\begin{align*}
    \disjunctioni{1} &:= \lprelaxation \cap \{\mathbf{x} \in \reals^{n} | \splitcoeff^{\intercal} \mathbf{x} \leq \splitrhs \} \\
    \disjunctioni{2} &:= \lprelaxation \cap \{\mathbf{x} \in \reals^{n} | \splitcoeff^{\intercal} \mathbf{x} \geq \splitrhs + 1 \}
\end{align*}
The disjunction is valid as $\splitcoeff^{\intercal} \xv$ must always take an integer value in a feasible solution to \eqref{eq:mip} due to the design of \splitcoeff.
We observe that the disjunctive set \disjunction can be written as the complement of a set \sset intersected with \lprelaxation, where \sset is defined as:
\begin{align}
    \sset := \{\mathbf{x} \in \reals^{n} \,|\, \splitrhs < \splitcoeff^{\intercal}\mathbf{x} < \splitrhs + 1\} \label{eq:sset}
\end{align}

\begin{figure}[h]
\centering
\begin{subfigure}[b]{0.475\textwidth}
    \centering
    \begin{tikzpicture}[scale=1]
    % Create the grid points
    \foreach \x/\y in {0/0, 0/1, 0/2, 0/3, 1/0, 1/1, 1/2, 1/3, 2/0, 2/1, 2/2, 2/3, 3/0, 3/1, 3/2, 3/3}
    {
    \fill[opacity=0.4] (\x,\y) circle (2pt);
    }
    \foreach \x/\y in {0/0, 0/1, 0/2, 1/1, 1/2, 2/1, 2/2}
    {
    \fill[opacity=1] (\x,\y) circle (2pt);
    }
    % Create the nodes of the polygon
    \node (a) at (2.5,1.5) {};
    \node (b) at (1.5,0.5) {};
    \node (c) at (0,0) {};
    \node (d) at (0,2) {};
    \node (e) at (0,2.8) {};
    \node (f) at (0.2,3) {};
    \node (g) at (1.5,2.5) {};
    
    \node (d11) at (2, 1) {};
    \node (d12) at (0, 1) {};
    \node (d21) at (2, 2) {};
    \node (d22) at (0, 2) {};
    
    % draw the big polygon    
    \draw[thick] (a.center) -- (b.center) -- (c.center) -- (d.center)  -- (e.center) -- (f.center) -- (g.center) -- cycle;
    % Draw the smaller disjunctions
    \filldraw[thick,fill=blue!60,fill opacity=0.4] (d11.center) -- (b.center) -- (c.center) -- (d12.center) -- cycle;
    \filldraw[thick,fill=blue!60,fill opacity=0.4] (d21.center) -- (d22.center) -- (d.center) -- (e.center) -- (f.center) -- (g.center) -- cycle;
    % Label the LP optimal point
    \fill[draw=black,fill=red] (a.center) circle (2pt);
    
    % Draw the S free set
    \draw[thick,dashed] (-0.5,1) -- (3.5,1);
    \draw[thick,dashed] (-0.5,2) -- (3.5,2);
    \node[label=above:{\sset}] (f) at (1.2,1.1) {};
    
    \node[label=above:{\disjunctioni{1}}] (g) at (0.8,0.2) {};
    \node[label=above:{\disjunctioni{2}}] (h) at (0.8,2) {};
    
    % Dummy nodes
    \node (xxx) at (2,3.5) {};
    \node (xxy) at (2,-0.5) {};
    
    \end{tikzpicture}
\label{fig:split_cut_1}
\end{subfigure}
\hfill
\begin{subfigure}[b]{0.475\textwidth}
\centering
    \begin{tikzpicture}[scale=1]
    \centering
    % Create the grid points
    \foreach \x/\y in {0/0, 0/1, 0/2, 0/3, 1/0, 1/1, 1/2, 1/3, 2/0, 2/1, 2/2, 2/3, 3/0, 3/1, 3/2, 3/3}
    {
    \fill[opacity=0.4] (\x,\y) circle (2pt);
    }
    \foreach \x/\y in {0/0, 0/1, 0/2, 1/1, 1/2, 2/1, 2/2}
    {
    \fill[opacity=1] (\x,\y) circle (2pt);
    }
    % Create the nodes of the polygon
    \node (a) at (2.5,1.5) {};
    \node (b) at (1.5,0.5) {};
    \node (c) at (0,0) {};
    \node (d) at (0,2) {};
    \node (e) at (0,2.8) {};
    \node (f) at (0.2,3) {};
    \node (g) at (1.5,2.5) {};
    
    \node (d11) at (2, 1) {};
    \node (d12) at (0, 1) {};
    \node (d21) at (2, 2) {};
    \node (d22) at (0, 2) {};
    
    % draw the big polygon    
    \draw[thick] (a.center) -- (b.center) -- (c.center) -- (d.center)  -- (e.center) -- (f.center) -- (g.center) -- cycle;
    % Draw the smaller disjunctions
    \filldraw[thick,fill=blue!60,fill opacity=0.4] (d11.center) -- (b.center) -- (c.center) -- (d12.center) -- cycle;
    \filldraw[thick,fill=blue!60,fill opacity=0.4] (d21.center) -- (d22.center) -- (d.center) -- (e.center) -- (f.center) -- (g.center) -- cycle;
    % Label the LP optimal point
    \fill[draw=black,fill=red] (a.center) circle (2pt);
    
    % Draw the S free set
    \draw[thick,dashed] (-0.5,1) -- (3.5,1);
    \draw[thick,dashed] (-0.5,2) -- (3.5,2);
    
    % Draw the cut
    \draw[ultra thick, red] (2, -0.5) -- (2, 3.5);
    % Dummy nodes
    \node (xxx) at (2,3.5) {};
    \node (xxy) at (2,-0.5) {};
    \end{tikzpicture}
\label{fig:split_cut_2}
\end{subfigure}
\caption{(Left) An example (simple) split. (Right) An example (simple) split cut.}
\label{fig:split_cut}
\end{figure}
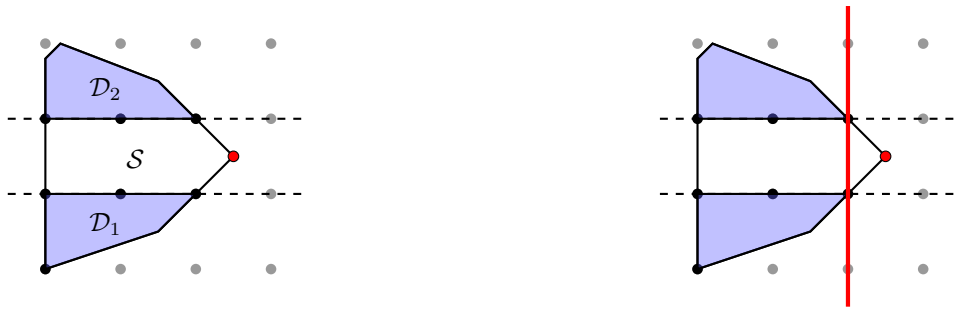

Note that notation is often abused, and the split can reference either the set \sset from \eqref{eq:sset} or the boundary of the set \sset, i.e.~the two hyperplanes from \eqref{eq:split}.
From a split disjunction, we can derive a \emph{split cut}. A split cut, \cut, is a valid inequality for both \disjunctioni{1} and \disjunctioni{2}, and separates some points from $\sset \cap \lprelaxation$. In the mixed-integer case, unlike the pure integer case \cite{schrijver1980cutting}, a finite amount of split cuts is not always sufficient for defining the integer hull and proving optimality, see \cite{cook1990chvatal}. A split is called \emph{simple} or \emph{elementary} if it only acts on a single variable, i.e.~$\splitcoeff = \mathbf{e}_{i}$ for some $i \in \mathcal{J}$. An example (simple) split disjunction alongside a valid split cut is visualised in Figure~\ref{fig:split_cut}.

\subsubsection{Intersection Cuts}
Some cuts are derived from the standard form of a MILP, which is defined using equality constraints instead of inequalities.
In particular, intersection cuts and GMI cuts are derived using this standard form.
Given our definition of a MILP in \eqref{eq:mip}, we can transform it to a standard form MILP in higher dimension by adding non-negative slack variables to each constraint. We can additionally substitute and introduce variables to shift variable bounds while keeping an equivalent formulation. We do this procedure to obtain the following MILP, where for ease of notation we will continue to use $\mathbf{c}$ and $\mathbf{A}$.
\begin{align}
    \underset{\mathbf{x}}{\text{argmin}}\{\mathbf{c}^{\intercal}\mathbf{x} \;\; | \;\; \mathbf{A}\mathbf{x} = \mathbf{b}, \;\; \mathbf{x} \geq \mathbf{0}, \;\; \mathbf{x} \in \mathbb{Z}^{|\mathcal{J}|} \times \mathbb{R}^{n + m - |\mathcal{J}|} \} \label{eq:mip_equality}
\end{align}

\begin{figure}[h]
\centering
\begin{subfigure}[b]{0.475\textwidth}
    \centering
    \begin{tikzpicture}[scale=1]
    % Create the grid points
    \foreach \x/\y in {0/0, 0/1, 0/2, 0/3, 1/0, 1/1, 1/2, 1/3, 2/0, 2/1, 2/2, 2/3, 3/0, 3/1, 3/2, 3/3}
    {
    \fill[opacity=0.4] (\x,\y) circle (2pt);
    }
    \foreach \x/\y in {1/0, 2/0, 2/1, 3/1, 3/2}
    {
    \fill[opacity=1] (\x,\y) circle (2pt);
    }
    % Create the nodes of the polygon
    \node (a) at (1.5,1.5) {};
    \node (b) at (3,2.5) {};
    \node (c) at (3,1) {};
    \node (d) at (2.5,0) {};
    \node (e) at (1.0,0) {};
    \node (f) at (1.0,0.75) {};
    
    % nodes for the cone
    \node (c1) at (3.5,2.83333) {};
    \node (c2) at (3.5,-0.5) {};
    \node (c3) at (0.166666,-0.5) {};
    
    % draw the big polygon    
    \draw[thick,fill=blue!60,fill opacity=0.4] (a.center) -- (b.center) -- (c.center) -- (d.center)  -- (e.center) -- (f.center) -- cycle;
    % Draw the cone
    \fill[pattern = north west lines] (a.center) -- (c1.center) -- (c2.center) -- (c3.center) -- cycle;
    % Label the LP optimal point
    \fill[draw=black,fill=red] (a.center) circle (2pt);
    % Draw the rays
    \draw[ultra thick,->] (a.center) -- (c1.center);
    \draw[ultra thick,->] (a.center) -- (c3.center);
    \node[label=above:{$\mathbf{r}_{1}$}] (g) at (0.5,0.3) {};
    \node[label=above:{$\mathbf{r}_{2}$}] (h) at (2.5,2.1) {};
    
    \end{tikzpicture}
\label{fig:lpcone_1}
\end{subfigure}
\hfill
\begin{subfigure}[b]{0.475\textwidth}
\centering
    \begin{tikzpicture}[scale=1]
    % Create the grid points
    \foreach \x/\y in {0/0, 0/1, 0/2, 0/3, 1/0, 1/1, 1/2, 1/3, 2/0, 2/1, 2/2, 2/3, 3/0, 3/1, 3/2, 3/3}
    {
    \fill[opacity=0.4] (\x,\y) circle (2pt);
    }
    \foreach \x/\y in {1/0, 2/0, 2/1, 3/1, 3/2}
    {
    \fill[opacity=1] (\x,\y) circle (2pt);
    }
    % Create the nodes of the polygon
    \node (a) at (1.5,1.5) {};
    \node (b) at (3,2.5) {};
    \node (c) at (3,1) {};
    \node (d) at (2.5,0) {};
    \node (e) at (1.0,0) {};
    \node (f) at (1.0,0.75) {};
    
    % nodes for the cone
    \node (coneapex) at (1.0,1.166666) {};
    \node (c1) at (3.5,2.83333) {};
    \node (c2) at (3.5,-0.5) {};
    \node (c3) at (1.0,-0.5) {};
    
    % draw the big polygon    
    \draw[thick,fill=blue!60,fill opacity=0.4] (a.center) -- (b.center) -- (c.center) -- (d.center)  -- (e.center) -- (f.center) -- cycle;
    % Draw the cone
    \fill[pattern = north west lines] (coneapex.center) -- (c1.center) -- (c2.center) -- (c3.center) -- cycle;
    % Label the LP optimal point
    \fill[draw=black,fill=red] (coneapex.center) circle (2pt);
    % Draw the rays
    \draw[ultra thick,->] (coneapex.center) -- (c1.center);
    \draw[ultra thick,->] (coneapex.center) -- (c3.center);
    \node[label=above:{$\mathbf{r}_{1}$}] (g) at (0.5,0.3) {};
    \node[label=above:{$\mathbf{r}_{2}$}] (h) at (2.5,2.1) {};
    
    \end{tikzpicture}
\label{fig:lpcone_2}
\end{subfigure}
\caption{The shaded area is the feasible region of \lpcone{\basissol}. The red dot is the apex of the simplicial conic relaxation of the feasible region of \eqref{eq:mip}, and $\mathbf{r}_{1}, \mathbf{r}_{2}$ are the rays of the cone. (Left) The red dot is both \lpoptimal and \basissol. (Right) The red dot is a primal infeasible \basissol.}
\label{fig:lpcone}
\end{figure}
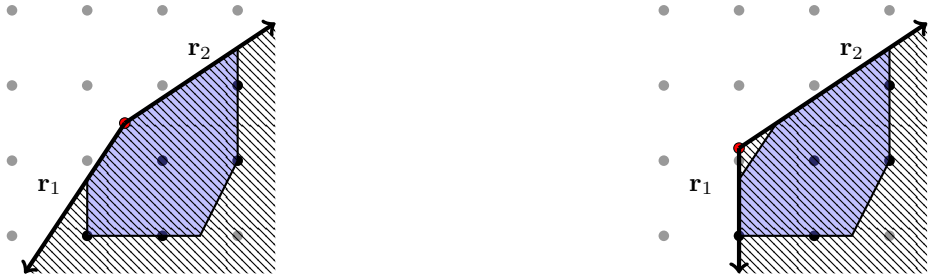

The simplex method typically used to solve LP relaxations of \eqref{eq:mip_equality} returns a \emph{basis}, $\mathcal{B} \subseteq \{1,...,n+m\}$, where $|\mathcal{B}| = m$. The basis is an index set of variables and relates to an extreme point of the LP relaxation, $\basissol \in \reals^{n+m}$, which is a basic solution. In practice the simplex method returns the optimal basic solution \lpoptimal. Associated with every basic solution \basissol is the LP cone, or corner polyhedron, $\lpcone{\basissol} \subseteq \reals^{n+m}$, whose apex is \basissol and whose rays are defined by the n-hyperplanes that form the basis. These rays are the columns of the simplex tableau relating to the non-basic variables. Note that in the case of primal degeneracy, multiple bases may result in the same extreme point but in different LP cones, and as such $\lpcone{\mathcal{B}}$ is the more appropriate notation. Two example LP cones are visualised in Figure~\ref{fig:lpcone}.

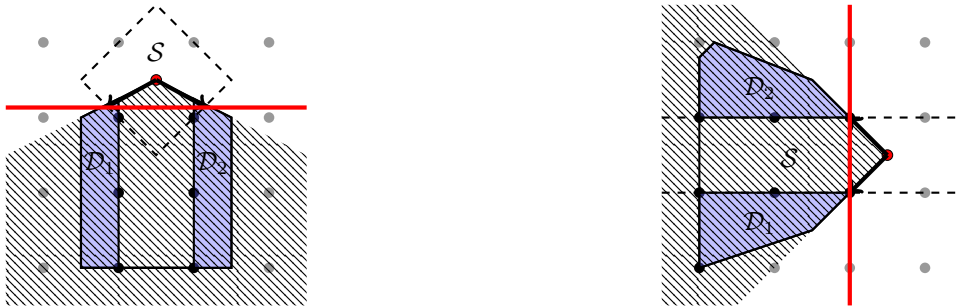
\begin{figure}[h]
\centering
\begin{subfigure}[b]{0.475\textwidth}
    \centering
    \begin{tikzpicture}[scale=1]
    % Create the grid points
    \foreach \x/\y in {0/0, 0/1, 0/2, 0/3, 1/0, 1/1, 1/2, 1/3, 2/0, 2/1, 2/2, 2/3, 3/0, 3/1, 3/2, 3/3}
    {
    \fill[opacity=0.4] (\x,\y) circle (2pt);
    }
    \foreach \x/\y in {1/0, 1/1, 1/2, 2/0, 2/1, 2/2}
    {
    \fill[opacity=1] (\x,\y) circle (2pt);
    }
    % Create the nodes of the polygon
    \node (a) at (2.5,2) {};
    \node (b) at (2.5,0) {};
    \node (c) at (0.5,0) {};
    \node (d) at (0.5,2) {};
    \node (e) at (1.5,2.5) {};
    
    \node (d11) at (1,2.25) {};
    \node (d12) at (1,0) {};
    \node (d21) at (2,2.25) {};
    \node (d22) at (2,0) {};
    
    % nodes for the cone
    \node (c1) at (3.5,1.5) {};
    \node (c2) at (3.5,-0.5) {};
    \node (c3) at (-0.5,-0.5) {};
    \node (c4) at (-0.5,1.5) {};
    
    % draw the big polygon    
    \draw[thick] (a.center) -- (b.center) -- (c.center) -- (d.center)  -- (e.center) -- cycle;
    % Draw the smaller disjunctions
    \filldraw[thick,fill=blue!60,fill opacity=0.4] (d21.center) -- (a.center) -- (b.center) -- (c.center) -- (d22.center) -- cycle;
    \filldraw[thick,fill=blue!60,fill opacity=0.4] (d11.center) -- (d12.center) -- (c.center) -- (d.center) -- cycle;
    % Label the LP optimal point
    \fill[draw=black,fill=red] (e.center) circle (2pt);
    
    % Draw the rays
    \draw[ultra thick,->] (e.center) -- (2.16666,2.13333);
    \draw[ultra thick,->] (e.center) -- (0.83333,2.13333);
    
    % Draw the cone
    \fill[pattern = north west lines] (e.center) -- (c1.center) -- (c2.center) -- (c3.center) -- (c4.center) -- cycle;
    
    % Draw the S free set
    \draw[thick,dashed] (0.5,2.5) -- (1.5,3.5) -- (2.5,2.5) -- (1.5,1.5) -- (0.5,2.5);
    \node[label=above:{\sset}] (f) at (1.5,2.5) {};
    
    \node[label=above:{\disjunctioni{1}}] (g) at (0.75,1) {};
    \node[label=above:{\disjunctioni{2}}] (h) at (2.25,1) {};
    
    % Draw the cut
    \draw[ultra thick, red] (-0.5, 2.13333) -- (3.5, 2.133333);
    
    \end{tikzpicture}
\label{fig:intersection_cut_1}
\end{subfigure}
\hfill
\begin{subfigure}[b]{0.475\textwidth}
\centering
    \begin{tikzpicture}[scale=1]
    \centering
    % Create the grid points
    \foreach \x/\y in {0/0, 0/1, 0/2, 0/3, 1/0, 1/1, 1/2, 1/3, 2/0, 2/1, 2/2, 2/3, 3/0, 3/1, 3/2, 3/3}
    {
    \fill[opacity=0.4] (\x,\y) circle (2pt);
    }
    \foreach \x/\y in {0/0, 0/1, 0/2, 1/1, 1/2, 2/1, 2/2}
    {
    \fill[opacity=1] (\x,\y) circle (2pt);
    }
    % Create the nodes of the polygon
    \node (a) at (2.5,1.5) {};
    \node (b) at (1.5,0.5) {};
    \node (c) at (0,0) {};
    \node (d) at (0,2) {};
    \node (e) at (0,2.8) {};
    \node (f) at (0.2,3) {};
    \node (g) at (1.5,2.5) {};
    
    \node (d11) at (2, 1) {};
    \node (d12) at (0, 1) {};
    \node (d21) at (2, 2) {};
    \node (d22) at (0, 2) {};
    
    % nodes for the cone
    \node (c1) at (0.5,-0.5) {};
    \node (c2) at (-0.5,-0.5) {};
    \node (c3) at (-0.5,3.5) {};
    \node (c4) at (0.5,3.5) {};
    
    % draw the big polygon    
    \draw[thick] (a.center) -- (b.center) -- (c.center) -- (d.center)  -- (e.center) -- (f.center) -- (g.center) -- cycle;
    % Draw the smaller disjunctions
    \filldraw[thick,fill=blue!60,fill opacity=0.4] (d11.center) -- (b.center) -- (c.center) -- (d12.center) -- cycle;
    \filldraw[thick,fill=blue!60,fill opacity=0.4] (d21.center) -- (d22.center) -- (d.center) -- (e.center) -- (f.center) -- (g.center) -- cycle;
    % Label the LP optimal point
    \fill[draw=black,fill=red] (a.center) circle (2pt);
    
    % Draw the rays
    \draw[ultra thick,->] (a.center) -- (2,2);
    \draw[ultra thick,->] (a.center) -- (2,1);
    
    % Draw the cone
    \fill[pattern = north west lines] (a.center) -- (c1.center) -- (c2.center) -- (c3.center) -- (c4.center) -- cycle;
    
    % Draw the S free set
    \draw[thick,dashed] (-0.5,1) -- (3.5,1);
    \draw[thick,dashed] (-0.5,2) -- (3.5,2);
    \node[label=above:{\sset}] (f) at (1.2,1.1) {};
    
    \node[label=above:{\disjunctioni{1}}] (g) at (0.8,0.2) {};
    \node[label=above:{\disjunctioni{2}}] (h) at (0.8,2) {};
    
    % Draw the cut
    \draw[ultra thick, red] (2, -0.5) -- (2, 3.5);
    % Dummy nodes
    \node (xxx) at (2,3.5) {};
    \node (xxy) at (2,-0.5) {};
    \end{tikzpicture}
\label{fig:intersection_cut_2}
\end{subfigure}
\caption{(Left) Example intersection cut. (Right) Example intersection cut that is also a split cut.}
\label{fig:intersection_cut}
\end{figure}

An intersection cut, similar to a split cut, is defined w.r.t.~a set $\sset \subseteq \reals^{n+m}$, which lies in the same dimension as $\mathbf{x}$ in the new space. Unlike split cuts however, \sset is not necessarily defined by two hyperplanes. Rather, it needs to be convex, to contain in its interior a current LP-feasible fractional solution we want to separate, and to not contain any integer-feasible solution in its interior.
In the context of MILP, the set \sset is a lattice-free set \cite{conforti2011corner}. In addition to the set \sset, an intersection cut is also defined w.r.t.~a simplicial conic relaxation of the feasible region of \eqref{eq:mip}, see Figure~\ref{fig:lpcone} for example simplicial conic relaxations derived from bases.
We note that while any simplicial conic relaxation can be exploited, an LP cone derived from a basis is used in practice. The idea to generate intersection cuts is to collect the intersection points of each ray with the boundary of the closure of \sset, and form a valid inequality as the hyperplane that contains all the intersection points. When using the LP cone \lpcone{\lpoptimal} and a set \sset containing \lpoptimal, the generated inequality will be a cut. Two examples of intersection cuts are visualised in Figure~\ref{fig:intersection_cut}.
For a deeper look into intersection cuts, we refer readers to \cite{conforti2011corner,balas2018disjunctive}.

\subsection{GMI Inequality Derivation}
\label{subsec:gmi_derive}

We will now derive the GMI inequality, which we note again is general-purpose and requires no additional problem structure. 

\begin{definition}[GMI inequality]
Given a valid equality for the LP relaxation of \eqref{eq:mip_equality}, $\mathbf{a}^{\intercal} \xv = b$, we distinguish the variables into those with integer requirements and those that are continuous, i.e., $\sum_{i \in \mathcal{J}} a_{i} x_{i} + \sum_{i \in [n] \setminus \mathcal{J}} a_{i}x_{i} = b$. Let $[n] = \{1,...,n\}$, $b = \floor{b} + f_{0}$, where $0 < f_{0} < 1$, and $a_{i} = \floor{a_{i}} + f_{i}$, where $0 \leq f_{i} < 1$ and $i \in [n]$. The GMI inequality is:
\begin{align}
    \sum_{i \in \mathcal{J}, f_{i} \leq f_{0}} \frac{f_{i}}{f_{0}}x_{i} + \sum_{i \in \mathcal{J}, f_{i} > f_{0}} \frac{1 - f_{i}}{1 - f_{0}}x_{i} + \sum_{i \in [n] \setminus \mathcal{J}, a_{i} \geq 0} \frac{a_{i}}{f_{0}}x_{i} - \sum_{i \in [n] \setminus \mathcal{J}, a_{i} < 0} \frac{a_{i}}{1 - f_{0}}x_{i} \geq 1 \label{eq:gomory}
\end{align}
\end{definition}
\begin{proof}[Derivation]
The logic of the GMI inequality is that if $f_{0} > 0$, then fractional multiples of integer variables and multiples of continuous variables must account for $f_{0}$. Specifically, they must sum to $f_{0}$ and a potential integer. That is:
\begin{align}
    \sum_{i \in \mathcal{J}, f_{i} \leq f_{0}} f_{i}x_{i} + \sum_{i \in \mathcal{J}, f_{i} > f_{0}} (f_{i} - 1)x_{i} + \sum_{i \in [n] \setminus \mathcal{J}} a_{i}x_{i} = k + f_{0}, \quad k \in \ints \label{eq:frac_sum}
\end{align}
This partition of $f_{i}$ values around $f_{0}$ is possible due to the observation that for $i$ such that $f_{i} > 0$, $a_i$ can be written as $a_{i} = \floor{a_{i}} + f_{i}$ or as $a_{i} = \ciel{a_{i}} + (f_{i} - 1)$. For example, $3.6 = 3 + 0.6$ or equivalently, $3.6 = 4 - 0.4$. This partition is done as it results in a strictly stronger cut than if the case $a_{i} = \floor{a_{i}} + f_{i}$ is always used \cite{cornuejols2008valid,andersen2005reduce,cornuejols2011improved}. Specifically, it results in smaller coefficients for like terms as $\frac{1-f_{i}}{1-f_{0}} < \frac{f_{i}}{f_{0}}$ when $f_{i} > f_{0}$.

Let us create a disjunction for two cases for inequality \eqref{eq:frac_sum},
where $k \leq -1$ or $k \geq 0$. In the case $k \leq -1$ we have that:
\begin{align}
    \sum_{i \in \mathcal{J}, f_{i} \leq f_{0}} f_{i}x_{i} + \sum_{i \in \mathcal{J}, f_{i} > f_{0}} (f_{i} - 1)x_{i} + \sum_{i \in [n] \setminus \mathcal{J}} a_{i}x_{i} &\leq -(1 - f_{0}) \nonumber \\
    \Rightarrow - \sum_{i \in \mathcal{J}, f_{i} \leq f_{0}} \frac{f_{i}}{1 - f_{0}}x_{i} + \sum_{i \in \mathcal{J}, f_{i} > f_{0}} \frac{1 - f_{i}}{1 - f_{0}}x_{i} - \sum_{i \in [n] \setminus \mathcal{J}} \frac{a_{i}}{1 - f_{0}}x_{i} &\geq 1 \label{eq:gomory_disjunction_1}
\end{align}
In the second case, $k \geq 0$ we have that:
\begin{align}
    \sum_{i \in \mathcal{J}, f_{i} \leq f_{0}} f_{i}x_{i} + \sum_{i \in \mathcal{J}, f_{i} > f_{0}} (f_{i} - 1)x_{i} + \sum_{i \in [n] \setminus \mathcal{J}} a_{i}x_{i} &\geq f_{0} \nonumber \\
    \Rightarrow \sum_{i \in \mathcal{J}, f_{i} \leq f_{0}} \frac{f_{i}}{f_{0}}x_{i} - \sum_{i \in \mathcal{J}, f_{i} > f_{0}} \frac{1 - f_{i}}{f_{0}}x_{i} + \sum_{i \in [n] \setminus \mathcal{J}} \frac{a_{i}}{f_{0}} x_{i} &\geq 1 \label{eq:gomory_disjunction_2}
\end{align}

As $\mathbf{x} \geq 0$, we can derive a globally valid inequality for the disjunctive set from the inequalities \eqref{eq:gomory_disjunction_1} - \eqref{eq:gomory_disjunction_2} by taking the maximum coefficient of each term over the two inequalities. That is the inequalities $\mathbf{a}^{\intercal}\xv \geq 1$ and $\mathbf{a}'^{\intercal}\xv \geq 1$ implies $\sum_{i=1}^{n}\max(a_{i},a'_{i})x_{i} \geq 1$.
We have grouped the terms in their derivation above such that at most one is positive.
The result of this derivation is exactly the GMI inequality \eqref{eq:gomory}.
\end{proof}

\subsection{GMI Cuts in Practice}
\label{subsec:gmi_practice}

In general, it is $\mathcal{NP}$-hard to find a GMI cut that separates a given LP-\emph{feasible} solution, or to determine if such a cut exists, see \cite{caprara2003separation,cornuejols2008valid}.
It is not $\mathcal{NP}$-hard, however, to separate a given basic solution of \lprelaxation, e.g., an LP-\emph{optimal} solution found by a simplex algorithm.

Consider a row of the simplex tableau for variable $x_{j}$ of basis $\mathcal{B}$. The row is an aggregated equality constraint, created from a linear combination of original constraints, where basic variable $x_{j}$ is described purely in terms of the non-basic variables. That is:
\begin{align}
    x_{j} = \basissoli{j} - \sum_{i \notin \mathcal{B}} \bar{a}_{ji} x_{i} \label{eq:simplex_row}
\end{align}
Here $\basissoli{j}$ is the right-hand side value of the tableau row and $\bar{a}_{ji}$ is the tableau entry for the row of basic variable $x_{j}$ and column of variable $x_{i}$. In the considered case of $\xv \geq \mathbf{0}$, the $\basissoli{j}$ is the value of variable $x_{j}$ at the basic solution, with the non-basic variables all taking values 0.

A GMI cut is derived from applying the GMI inequality procedure from Subsection~\ref{subsec:gmi_derive} to the aggregated equality constraint \eqref{eq:simplex_row}. This procedure is only applied to rows of the simplex tableau that correspond to integer variables with fractional LP solutions. This is because these rows have a fractional right hand side, and the resulting GMI inequality guarantees separation of the current LP solution. An inequality produced by this method is called a GMI cut. 

Geometrically, a GMI cut is a split cut, and therefore also an intersection cut and a disjunctive cut.
Specifically, it is an intersection cut for the split $\disjunction(\splitcoeff^{G},\floor{\basissoli{j}})$, where $\splitcoeff^{G}$ is defined as follows:
\begin{align}
    \splitcoeff^{G} \in \ints^{n}, \quad \text{where } \splitcoeff^{G}_{i} &:= 
    \begin{cases}
        \floor{\bar{a}_{ji}}, &\quad \text{if } (f_{i} \leq f_{0}) \; \land \; i \notin \mathcal{B} \\
        \ciel{\bar{a}_{ji}}, &\quad \text{if } (f_{i} > f_{0}) \; \land \; i \notin \mathcal{B} \\
        1, &\quad \text{if } i = j \\
        0, &\quad \text{if } (i \neq j) \land i \in \mathcal{B}\\
    \end{cases}
    \quad \forall i \in \{1,...,n\}
\end{align}

The GMI cut of the tableau row \eqref{eq:simplex_row} is the strengthened version of the intersection cut obtained from the elementary split $\disjunction(\mathbf{e}_{j}, \floor{\basissoli{j}})$, see \cite{andersen2005reduce,balas2018disjunctive}. Deriving a cut using the elementary split for the simplex tableau row \eqref{eq:simplex_row} of variable $x_{j}$ without the strengthening procedure results in the intersection cut \eqref{eq:weak_gomory}. This cut is obtained by treating integer variables similarly to continuous ones for the GMI derivation.
\begin{align}
    \sum_{i \notin \mathcal{B}, \bar{a}_{ji} \geq 0} \frac{\bar{a}_{ji}}{f_{0}}x_{i} - \sum_{i \notin \mathcal{B}, \bar{a}_{ji} < 0} \frac{\bar{a}_{ji}}{1 - f_{0}}x_{i} \geq 1 \label{eq:weak_gomory}
\end{align}
We denote this inequality as \emph{weak-GMI}, and note that the GMI cut will always dominate the associated weak-GMI cut. We also note that the strengthening procedure is performed by using fractional coefficient values $f_{i}$ instead of $\bar{a}_{ij}$ for integer variables and the partitioning of those fractional coefficients $f_{i}$ around $f_{0}$.

\section{Cutting Plane Selection for Variable Selection}

The core idea of our work is to use measures of cuts to evaluate and decide on corresponding branching candidates. Specifically, we will generate the GMI cut from the corresponding tableau row of each branching candidate, and use cut selection techniques to dictate branching decisions. We will additionally augment the default SCIP hybrid branching rule with history-based scores of already-computed GMI cuts from previous separation rounds. 

Currently, history-based approaches, see \cite{achterberg2005branching,achterberg2009hybrid}, are the backbone behind branching rules used in MILP solvers \cite{scip8,highs}. \emph{Pseudo-costs} \cite{pseudocost}, the most prolific case of history-based approaches, estimate scores for a branching candidate based on the historical objective value improvement of child nodes spawned from branching on the candidate. One can consider pseudo-costs as an approximation of \emph{strong-branching} scores, see e.g.~\cite{achterberg2005branching}, which are derived from directly solving the upper and lower LP relaxations of all branching candidates. In our approach, branching scores are derived from cut quality measures of cuts generated from each branching candidate. It is complementary to pseudo-costs in that it provides an objective-free measure. These cut-based scores can be integrated into SCIP's default scoring rule, using a history of cut quality measures from previously generated cuts. This is similar to other history-based scores, such as those based on bound inferences, conflict information, and subproblem infeasibility~\cite{achterberg2009hybrid}. 

The classical cut scoring measure is \emph{efficacy}\footnote{Main selection criteria for most MILP solvers, e.g., FICO Xpress 9.0 and SCIP 8.0}, which denotes the Euclidean distance between the LP optimal solution and the cut hyperplane. Given a cut $(\coefficients, \beta) \in \reals^{n+1}$ and the LP optimal solution \lpoptimal, efficacy is defined as:
\begin{align}
    \eff{\coefficients}{\beta}{\lpoptimal} := \frac{\coefficients^{\intercal} \lpoptimal  - \beta}{\|\coefficients\| } \label{eq:efficacy}
\end{align}

When scoring cuts for the purpose of branching, we will rely on efficacy as our cut measure. We note that there exists many more potential cut scoring measures \cite{wesselmann2012implementing,turner2022cutting}, however preliminary results of their inclusion led to negligible improvements.

GMI cuts are not the only cuts associated with split disjunctions, or even the elementary split, i.e.~branching decisions. For example, lift-and-project cuts \cite{balas1979disjunctive, perregaard2003generating} are intersection cuts of elementary splits. The elementary splits from which these cuts are derived, however, are not necessarily related to the current LP basis, nor even necessarily related to a primal-feasible LP basis. Nevertheless, scoring measures for this family of cuts are also a potentially potent indicator of good branching decisions. We however restricted our study to GMI cuts which are readily computed and available for all variables in all MILP solvers.

\section{Experiments}
 
We conduct three experiments: first, we analyse the effectiveness of our initial approach compared to standard branching rules (Subsection \ref{subsec:experiment_branch_rules}). Then, we refine our approach to a history-based one, and determine the best parameter value for including our approach in the state-of-the-art branching rule \emph{hybrid branching} (Subsection \ref{subsec:gmi_hist_branching}). Finally,  we compare our integrated branching rule to default SCIP with experiments run in exclusive mode, i.e. one job per machine (Subsection \ref{subsec:improve_default}). We perform experiments on the MIPLIB 2017 benchmark set\footnote{MIPLIB 2017 -- The Mixed Integer Programming Library \url{https://miplib.zib.de/}.} \cite{miplib}, which we will now simply refer to as MIPLIB. For all these experiments we present two variants: Firstly, we use default SCIP on the original instances to analyse the impact on the the out-of-the-box behaviour of a MILP solver. Secondly, we use SCIP with heuristics disabled and the optimal solution provided, which reduces random noise and emphasises the effect of branching rules. 

We define a run as an instance random-seed pair for which we use a given branching rule. All results are obtained by averaging results over the SCIP random seeds $\{1,2,3,4,5\}$. For all experiments, SCIP 8.0.3 \cite{scip8} is used, with PySCIPOpt \cite{pyscipopt} as the API, and Xpress 9.0.2 \cite{xpress} as the LP solver. For Subsections \ref{subsec:experiment_branch_rules} and \ref{subsec:gmi_hist_branching}, experiments are run in non-exclusive mode on a cluster equipped with Intel Xeon Gold 6342 CPUs running at 2.80GHz, where each run is restricted to 2GB memory, 2h time limit, and the LP solver is restricted to a single thread. For Subsection \ref{subsec:improve_default} experiments are run in exclusive mode on a cluster equipped with Intel Xeon Gold 5122 CPUs running at 3.60GHz, where each run is restricted to 48GB memory, 2h time limit, and the LP solver is restricted to a single thread. The code used for all experiments is available and open-source\footnote{\url{https://github.com/Opt-Mucca/branching-via-cut-selection}}, and will be integrated in the next release of SCIP.

For the entirety of our experiments, we filter out any instance that for any random seed was solved to optimality without branching, hit a memory limit, or encountered LP errors. Note that the instance is only filtered in a comparison of branching rules when one of the criteria is met for a run on one of the compared branching rules. When comparing results from branching rules, we use individual instance-seed pairs as data points as opposed to the aggregate performance over the random seeds. Additionally, when shifted geometric means are referenced, we use a shift of 100, 10s, and 1s for number of nodes, solving time, and branching time respectively. We finally note that certain instances were excluded when no optimal solution was available on the MIPLIB website.

\subsection{Gomory Cut-Based Branching Rules}\label{subsec:experiment_branch_rules}

To rank the effectiveness of our GMI-based branching rules, we compare them against standard branching rules from the literature, with Table~\ref{tab:branch_rules} containing a complete list.

\begin{table}[h]
\centering
\begin{tabular}{lc}
\hline
Branching Rule & Description \\
\hline
\emph{GMI} & Generate GMI cuts from tableau. Select candidate from cut with largest efficacy. \\
\emph{weak-GMI} & Generate weak-GMI cuts from Tableau. Select candidate from cut with largest efficacy. \\
\emph{fullstrong} & Solve LP relaxations of children nodes for all candidates, see \cite{achterberg2005branching}. \\
\emph{hybrid} & Reliability pseudo-cost / Hybrid. (Default SCIP scoring rule, see \cite{achterberg2005branching,achterberg2009hybrid}) \\
\emph{random} & Select random candidate. \\
\hline
\end{tabular}
\caption{Branching rules used in Experiment \ref{subsec:experiment_branch_rules}}
\label{tab:branch_rules}
\end{table}

\begin{table}[h]
    \begin{subtable}{\textwidth}
    \centering
    \begin{tabular}{lc|ccccc}
        \hline
        \hline
        Metric & Pairs & \emph{GMI} & \emph{weak-GMI} & \emph{fullstrong} & \emph{hybrid} & \emph{random} \\
        \hline
        \hline
        \noalign{\vskip 1mm}
        \multicolumn{7}{c}{MIPLIB} \\
        \noalign{\vskip 1mm}
        \hline
        \noalign{\vskip 0.35mm}
        Nodes & 237 & 3877 & 3109 & \textbf{512} & 1562 & 8589 \\
        Time (s) & 237 & 222 & 191 & 289 & \textbf{101} & 214 \\
        Time w/o branch time (s) & 237 & 124 & 111 & \textbf{64} & 81 & 214 \\
        Branch time (s) & 237 & 47 & 39 & 107 & 11 & \textbf{0} \\
        \hline
        \noalign{\vskip 1mm}
        \multicolumn{7}{c}{MIPLIB with optimal solution provided and heuristics disabled} \\
        \noalign{\vskip 1mm}
        \hline
        \noalign{\vskip 0.35mm}
        Nodes & 227 & 3767 & 2812 & \textbf{386} & 1287 & 8286 \\
        Time (s) & 227 & 179 & 138 & 126 & \textbf{73} & 168 \\
        Time w/o branch time (s) & 227 & 97 & 82 & \textbf{44} & 59 & 168 \\
        Branch time (s) & 227 & 46 & 32 & 39 & 8 & \textbf{0} \\
        \hline
        \hline
    \end{tabular}
    \caption{Instance-seed pairs where all branching rules solved to optimality.}
    \end{subtable}
    \begin{subtable}{\textwidth}
    \centering
    \begin{tabular}{lc|ccccc}
        \hline
        \hline
        Metric & Pairs & \emph{GMI} & \emph{weak-GMI} & \emph{fullstrong} & \emph{hybrid} & \emph{random} \\
        \hline
        \hline
        \noalign{\vskip 1mm}
        \multicolumn{7}{c}{MIPLIB} \\
        \noalign{\vskip 1mm}
        \hline
        \noalign{\vskip 0.35mm}
        Time (s) & 464 & 977 & 874 & 1200 & \textbf{310} & 858 \\
        Time w/o branch time (s) & 464 & 362 & 337 & \textbf{113} & 260 & 858 \\
        Branch time (s) & 464 & 340 & 298 & 689 & 26 & \textbf{1} \\
        \hline
        \noalign{\vskip 1mm}
        \multicolumn{7}{c}{MIPLIB with optimal solution provided and heuristics disabled} \\
        \noalign{\vskip 1mm}
        \hline
        \noalign{\vskip 0.35mm}
        Time (s) & 440 & 1016 & 832 & 610 & \textbf{269} & 943 \\
        Time w/o branch time (s) & 440 & 338 & 294 & \textbf{85} & 225 & 942 \\
        Branch time (s) & 440 & 401 & 311 & 306 & 24 & \textbf{1} \\
        \hline
        \hline
    \end{tabular}
    \caption{Instance-seed pairs where at-least one branching rule solved to optimality.}
    \end{subtable}
    \caption{Shifted geometric mean results. Best branching rule per metric in \textbf{bold}.}
    \label{tab:gmi_geom_mean}
\end{table}

The shifted geometric means over three performance metrics on our data sets are presented in Table \ref{tab:gmi_geom_mean}. We observe expected performance from the standard branching rules. \emph{Fullstrong} requires the least nodes to prove optimality over all data sets, while \emph{hybrid} is the branching rule that most quickly proves optimality. Our newly introduced branching rule \emph{GMI}, is regrettably inferior to default SCIP over all metrics and data sets, however we observe that it clearly has a positive signal due to it requiring substantially less nodes than \emph{random} to prove optimality over all data sets. Most interesting is the relative performance of \emph{weak-GMI} to \emph{GMI}, where \emph{weak-GMI} wins over all metrics and data sets. This suggests that the strengthened cut, while strictly better than the weaker version in a cutting plane context, has lost some level of the representation of the disjunction that the weaker cut is derived from.

For running time, we must also address the overhead of our branching rule. While ultimately faster per node than strong branching, we still need to generate a GMI cut for every branching candidate at every node. This overhead is significant, and is the reason why \emph{random} is on average faster to solve over MIPLIB both with and without a provided solution. This is despite requiring over twice as many nodes. This can be verified by seeing that \emph{GMI} and \emph{weak-GMI} both are much faster than \emph{random} when removing branching time from consideration.

\subsection{History-Based GMI Branching} \label{subsec:gmi_hist_branching}

Our approaches \emph{GMI} and \emph{weak-GMI} were shown to make substantially better branching decisions than \emph{random}, but were ultimately too slow, and were not as good as LP relaxation-based branching rules. The default SCIP branching rule, while dominated by pseudo-costs, is a hybrid method, with scores from a weighted sum of metrics. Most of these metrics are history-based, meaning that they use information from different parts of the solving process, and are quick to evaluate. Given that GMI cuts are already generated by SCIP throughout the solve process, we store for each variable, the normalised efficacy of the most recent GMI cut generated from a tableau row when the variable is basic and fractional. This normalised efficacy can then be used to augment the branching candidate's score of default SCIP. We normalise efficacy by the maximum GMI cut's efficacy from the given separation round, and note that we only store the normalised efficacy of the most recent cut if the non-normalised efficacy is above some epsilon tolerance. We stress here that this approach requires no additional overhead, as the cuts themselves as well as their efficacies are already computed in the separation process. 

We denote our new branching rule \emph{gmi-}$10^{-x}$, where $10^{-x}$
denotes the coefficient used in the weighted sum scoring rule for the normalised efficacy of the last generated GMI. The shifted geometric mean of performance metrics for various coefficient values are presented in Table \ref{tab:hybrid_geom_mean}. We observe that too high of a coefficient, as in \emph{gmi-}$10^{-2}$, results in worse performance than default SCIP over all metrics and all data sets. By decreasing the coefficient value, we see an improvement in performance, with \emph{gmi-}$10^{-5}$ being the best performing rule w.r.t.~both nodes and solve time over all data sets. We also observe that the branching rules on either side of \emph{gmi-}$10^{-5}$, i.e.~\emph{gmi-}$10^{-4}$ and \emph{gmi-}$10^{-6}$, always outperform default SCIP, indicating that $10^{-5}$ is a sweet spot. We therefore conclude that $10^{-5}$ is a good and robust coefficient choice for improving hybrid branching, i.e.~default SCIP, once again noting that it requires no additional overhead since branching time is functionally identical.

\begin{table}[h]
    \begin{subtable}{\textwidth}
    \centering
    \begin{tabular}{lc|cccccc}
        \hline
        \hline
        Metric & Pairs & \emph{hybrid} & \emph{gmi-}$10^{-2}$ & \emph{gmi-}$10^{-3}$ & \emph{gmi-}$10^{-4}$ & \emph{gmi-}$10^{-5}$ & \emph{gmi-}$10^{-6}$ \\
        \hline
        \hline
        \noalign{\vskip 1mm}
        \multicolumn{8}{c}{MIPLIB} \\
        \noalign{\vskip 1mm}
        \hline
        \noalign{\vskip 0.35mm}
        Nodes & 472 & 4345 & 5025 & 4590 & \textit{4214} & \textbf{\textit{4101}} & \textit{4267} \\
        Time (s) & 472 & 253 & 276 & 253 & \textit{235} & \textbf{\textit{232}} & \textit{239} \\
        Time w/o branch time (s) & 472 & 208 & 228 & 209 & \textit{194} & \textbf{\textit{191}} & \textit{197} \\
        Branch time (s) & 472 & 24 & 26 & 24 & \textbf{\textit{22}} & \textit{23} & \textit{23} \\
        \hline
        \noalign{\vskip 1mm}
        \multicolumn{8}{c}{MIPLIB with optimal solution provided and heuristics disabled} \\
        \noalign{\vskip 1mm}
        \hline
        \noalign{\vskip 0.35mm}
        Nodes & 426 & 3870 & 4426 & 4023 & \textit{3822} & \textbf{\textit{3671}} & \textit{3820} \\
        Time (s) & 426 & 204 & 226 & 213 & \textit{203} & \textbf{\textit{195}} & \textit{199} \\
        Time w/o branch time (s) & 426 & 165 & 185 & 174 & 165 & \textbf{\textit{158}} & \textit{161} \\
        Branch time (s) & 426 & 22 & 23 & 22 & \textbf{\textit{21}} & \textbf{\textit{21}} & 22 \\
        \hline
        \hline
    \end{tabular}
    \caption{Instance-seed pairs where all branching rules solved to optimality.}
    \end{subtable}
    \begin{subtable}{\textwidth}
    \centering
    \begin{tabular}{lc|cccccc}
        \hline
        \hline
        Metric & Pairs & \emph{hybrid} & \emph{gmi-}$10^{-2}$ & \emph{gmi-}$10^{-3}$ & \emph{gmi-}$10^{-4}$ & \emph{gmi-}$10^{-5}$ & \emph{gmi-}$10^{-6}$ \\
        \hline
        \hline
        \noalign{\vskip 1mm}
        \multicolumn{8}{c}{MIPLIB} \\
        \noalign{\vskip 1mm}
        \hline
        \noalign{\vskip 0.35mm}
        Time (s) & 548 & 370 & 402 & 378 & \textit{350} & \textbf{\textit{342}} & \textit{348} \\
        Time w/o branch time (s) & 548 & 308 & 336 & 315 & \textit{293} & \textbf{\textit{284}} & \textit{290} \\
        Branch time (s) & 548 & 33 & 34 & \textit{32} & \textbf{\textit{31}} & \textbf{\textit{31}} & \textbf{\textit{31}} \\
        \hline
        \noalign{\vskip 1mm}
        \multicolumn{8}{c}{MIPLIB with optimal solution provided and heuristics disabled} \\
        \noalign{\vskip 1mm}
        \hline
        \noalign{\vskip 0.35mm}
        Time (s) & 507 & 310 & 357 & 333 & 314 & \textbf{\textit{293}} & \textit{301} \\
        Time w/o branch time (s) & 507 & 254 & 295 & 275 & 259 & \textbf{\textit{240}} & \textit{247} \\
        Branch time (s) & 507 & 31 & 34 & 32 & \textbf{\textit{30}} & \textbf{\textit{30}} & 31 \\
        \hline
        \hline
    \end{tabular}
    \caption{Instance-seed pairs where at-least one branching rule solved to optimality.}
    \end{subtable}
    \caption{Shifted geometric mean results. Branching rules better than default in \textit{italics}, best in \textbf{bold}.}
    \label{tab:hybrid_geom_mean}
\end{table}

Instead of simply using the efficacy of the most recently generated GMI cut, we performed many preliminary experiments using the historical normalised average efficacy of all generated GMI cuts.
%\footnote{In the original preprint we claimed to use the historical average. Due to a bug in the implementation we foolishly always actually retrieved the efficacy of the most recent GMI cut. (Edited out by Mathieu)}
While this approach also had parameter values that outperformed default SCIP, it was thoroughly outperformed by the efficacy of the most-recently generated GMI cut. For the implementation of both approaches, ignoring GMI cuts that only separated the LP solution by a marginal efficacy improved performance.
We also performed preliminary experiments using a new homogeneous MILP instance set, SNDlib-MIPs \cite{zenodo_sndlib}, which was inspired by SNDLib \cite{sndlib}, but that the results while better than default SCIP, were only a marginal improvement.

Up to this point, our runs were affected by our experimental setup, where memory limits reduced the size of instances that we considered, and the non-exclusive mode of the computation cluster introduced additional noise w.r.t.~solve time. We therefore perform a more in-depth comparison of \emph{hybrid} and \emph{gmi-}$10^{-5}$ over MIPLIB in the following subsection.

\subsection{Improving Default SCIP}\label{subsec:improve_default}

\begin{table}[h]
    \begin{subtable}{\textwidth}
    \centering
    \begin{tabular}{lc|cc}
        \hline
        \hline
        Metric & Pairs & \emph{hybrid} & \emph{gmi-}$10^{-5}$ \\
        \hline
        \hline
        \noalign{\vskip 1mm}
        \multicolumn{4}{c}{MIPLIB} \\
        \noalign{\vskip 1mm}
        \hline
        \noalign{\vskip 0.35mm}
        Nodes & 541 & 4385 & \textbf{4150} \\
        Time (s) & 541 & 349 & \textbf{339} \\
        \hline
        \noalign{\vskip 1mm}
        \multicolumn{4}{c}{MIPLIB with optimal solution provided and heuristics disabled} \\
        \noalign{\vskip 1mm}
        \hline
        \noalign{\vskip 0.35mm}
        Nodes & 491 & 4702 & \textbf{4334} \\
        Time (s) & 491 & 325 & \textbf{304} \\
        \hline
        \hline
    \end{tabular}
    \caption{Instance-seed pairs where all branching rules solved to optimality.}
    \end{subtable}
    \begin{subtable}{\textwidth}
    \centering
    \begin{tabular}{lc|cc}
        \hline
        \hline
        Metric & Pairs & \emph{hybrid} & \emph{gmi-}$10^{-5}$ \\
        \hline
        \hline
        \noalign{\vskip 1mm}
        \multicolumn{4}{c}{MIPLIB} \\
        \noalign{\vskip 1mm}
        \hline
        \noalign{\vskip 0.35mm}
        Time (s) & 578 & 411 & \textbf{398} \\
        \hline
        \noalign{\vskip 1mm}
        \multicolumn{4}{c}{MIPLIB with optimal solution provided and heuristics disabled} \\
        \noalign{\vskip 1mm}
        \hline
        \noalign{\vskip 0.35mm}
        Time (s) & 528 & 381 & \textbf{361} \\
        \hline
        \hline
    \end{tabular}
    \caption{Instance-seed pairs where at-least one branching rule solved to optimality.}
    \end{subtable}
    \caption{Shifted geometric mean results. Best branching rule per metric in \textbf{bold}.}
    \label{tab:head_to_head_geom_mean}
\end{table}

Our in-depth comparison presented in Table~\ref{tab:head_to_head_geom_mean} shows that our branching rule clearly outperforms default SCIP's hybrid branching. Over MIPLIB both with and without a solution provided, our augmented branching method results in faster solve times, and fewer nodes. We stress here, however, that due to the size of the coefficient value for the normalised efficacy of the last generated GMI cut, our approach will often act more as a tie breaker rather than as the predominant branching decision.

\begin{table}[h]
    \centering
    \begin{tabular}{cccccc}
        \hline
        \hline
        \noalign{\vskip 1mm}
        \multicolumn{6}{c}{MIPLIB} \\
        \multicolumn{6}{c}{67.1\% instance-seed pairs affected} \\
        \noalign{\vskip 1mm}
        \hline
        \noalign{\vskip 0.35mm}
        Time (s) & Nodes & $\Delta-$Solved & $\Delta-$Dual & $\Delta-$Primal \\
        \hline
        \textbf{0.96} & \textbf{0.92} & -1(/578) & -11(/472) & -4(/472) \\
        \hline
        \noalign{\vskip 1mm}
        \multicolumn{6}{c}{MIPLIB with optimal solution provided and heuristics disabled} \\
        \multicolumn{6}{c}{69.5\% instance-seed pairs affected} \\
        \noalign{\vskip 1mm}
        \hline
        \noalign{\vskip 0.35mm}
        \textbf{0.91} & \textbf{0.89} & -1(/528) & -6(/307) & - \\
        \hline
        \hline
    \end{tabular}
    \caption{Ratio of shifted geometric means for \emph{gmi-}$10^{-5}$ vs
    \emph{hybrid} over affected instances (solved to optimality for both branching rules). $\Delta$ is the difference in wins over the instance-seed pairs for amount solved, and the respective bounds for unsolved instances. Entries are in bold when our approach is better than SCIP default.}
    \label{tab:head_to_head_geom_ratio}
\end{table}

The performance improvement of \emph{gmi-}$10^{-5}$ becomes even more apparent when we look only at affected instances, see Table~\ref{tab:head_to_head_geom_ratio}. Over MIPLIB, we achieve a $4\%$ speedup on affected instances, and require $8\%$ fewer nodes. This improvement becomes even more apparent when an optimal solution is provided. Our approach however is outperformed by default on unsolved instances. We do stress that upon further investigation, we found no evidence that \emph{gmi-}$10^{-5}$ is outperformed on larger or more complex instances that solved within the time limit. The large amount of affected instances indicates that the MILP solver ends up in situations where the scores of branching candidates, especially the pseudo-costs, are very similar. 

\begin{table}[h]
    \centering
    \begin{tabular}{ccccc}
        \hline
        \hline
        \noalign{\vskip 1mm}
        \multicolumn{5}{c}{All instance-seed pairs} \\
        \noalign{\vskip 1mm}
        \hline
        \noalign{\vskip 0.35mm}
        pseudo-cost & conflict & inference & cutoff & gmi \\
        \hline
        0.208 (0.428) & 0.0 (0.084) & 1 (0.920) & 0 (0) & 0.544 (0.537) \\
        \hline
        \noalign{\vskip 1mm}
        \multicolumn{5}{c}{Affected instance-seed pairs (10.7\%)} \\
        \noalign{\vskip 1mm}
        \hline
        \noalign{\vskip 0.35mm}
        0.833 (0.562) & 0.0 (0.011) & 1 (0.873) & 0 (0) & 0.580 (0.556) \\
        \hline
        \hline 
    \end{tabular}
    \caption{Median (Mean) proportion of fractional variables with branching scores initialised at root node over MIPLIB with optimal solution provided and heuristics disabled.. A variable is marked as initialised if there is at least one record of the appropriate score. Affected instance-seed pairs are those with different root branching decisions for \emph{hybrid} and \emph{gmi-}$10^{-5}$.}
    \label{tab:changes_in_relpscost}
\end{table}

Given the frequency of affected instance-seed pairs when comparing our branching rule to SCIP default, we do an analysis of the initialisation of the different measures in hybrid branching at the root node. The results are presented in Table  \ref{tab:changes_in_relpscost}, where we reference \cite{achterberg2009hybrid} for definitions of conflict, inference, and cutoff. We observe that $10.7\%$ of instance-seed pairs make a different root node branching decision than \emph{hybrid} when using \emph{gmi-}$10^{-5}$ as a tiebreaker, and that $54.4\%$ of variables have at least one GMI cut generated from tableau rows when they were fractional and basic. This initialisation ratio is important, as it is substantially larger than that of pseudo-costs over the whole data set, and provides a new signal for branching that is well initialised before branching begins. We also note that affected instances have a much higher initialisation rate of pseudo-costs than average, indicating that a frequent issue is not necessarily the lack of strong branching initialisation, but rather that the best initialisations all take the same value. 

\begin{figure}[h]
    \centering
    \begin{subfigure}[b]{\linewidth}
    \centering
    \includegraphics[height=6.5cm,width=0.495\linewidth]{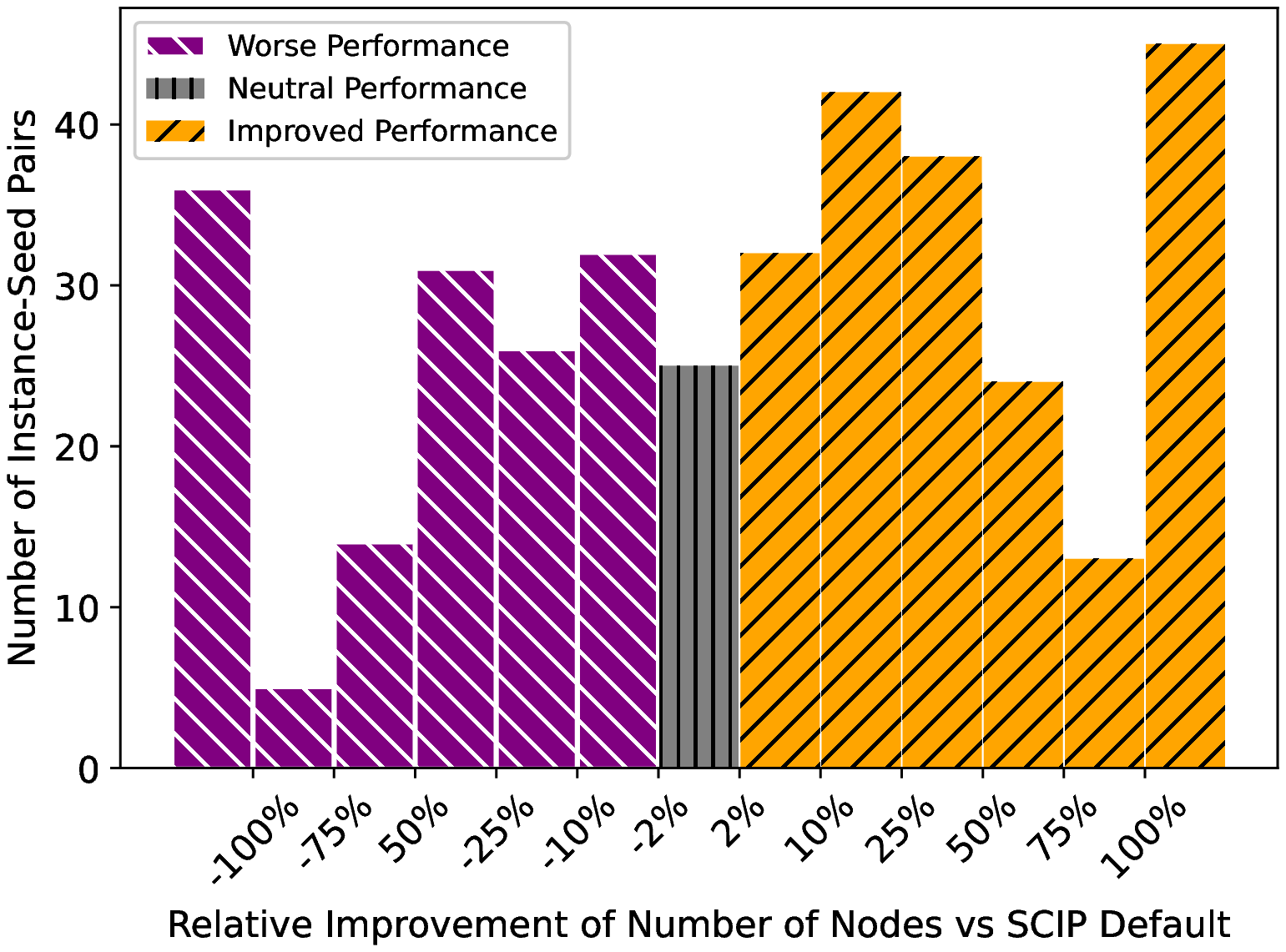}\hfill
    \includegraphics[height=6.5cm,width=0.495\linewidth]{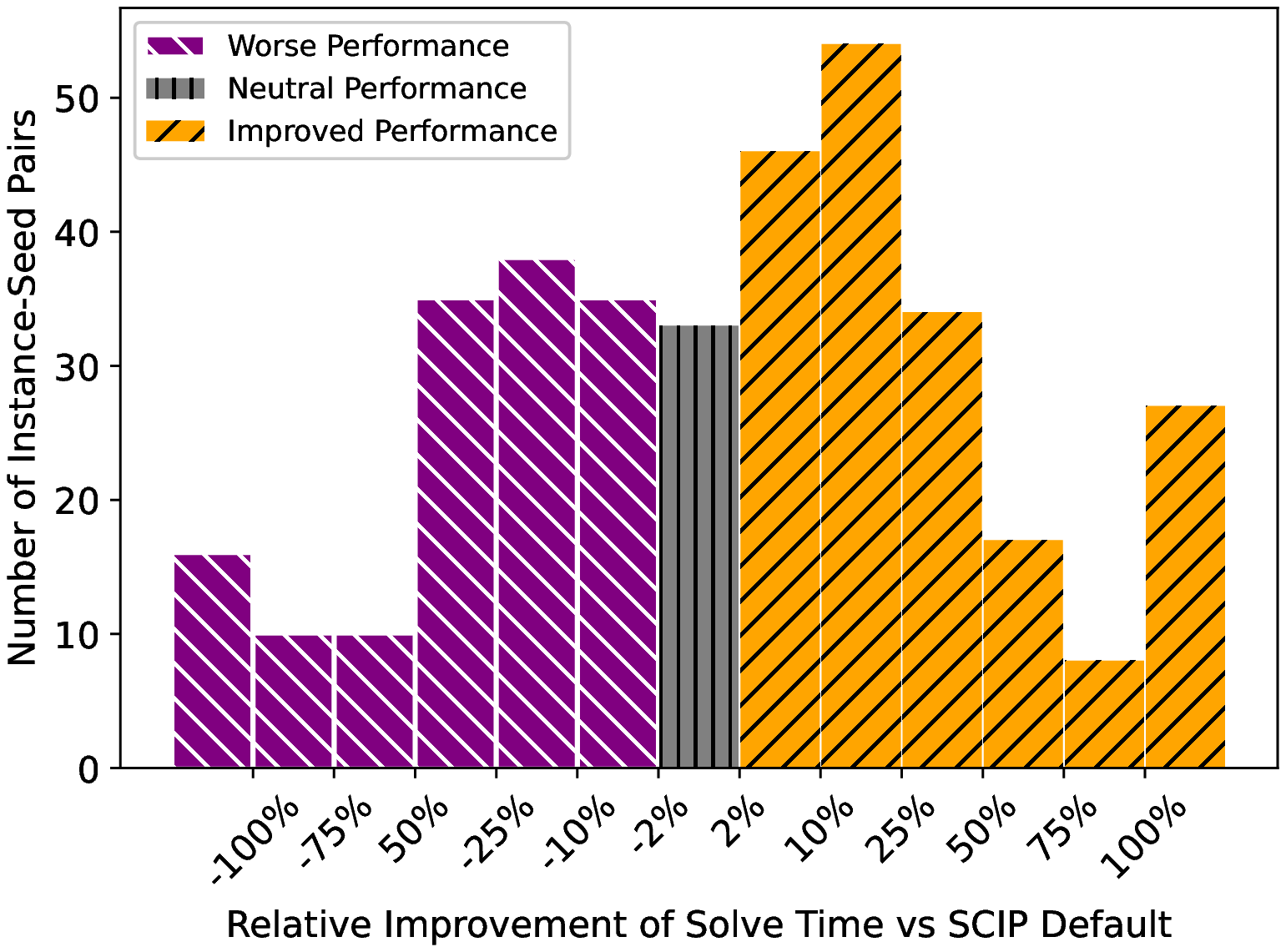}
    \caption{MIPLIB.}
    \end{subfigure}
    \vskip\baselineskip
    \begin{subfigure}[b]{\linewidth}
    \centering
    \includegraphics[height=6.5cm,width=0.48\linewidth]{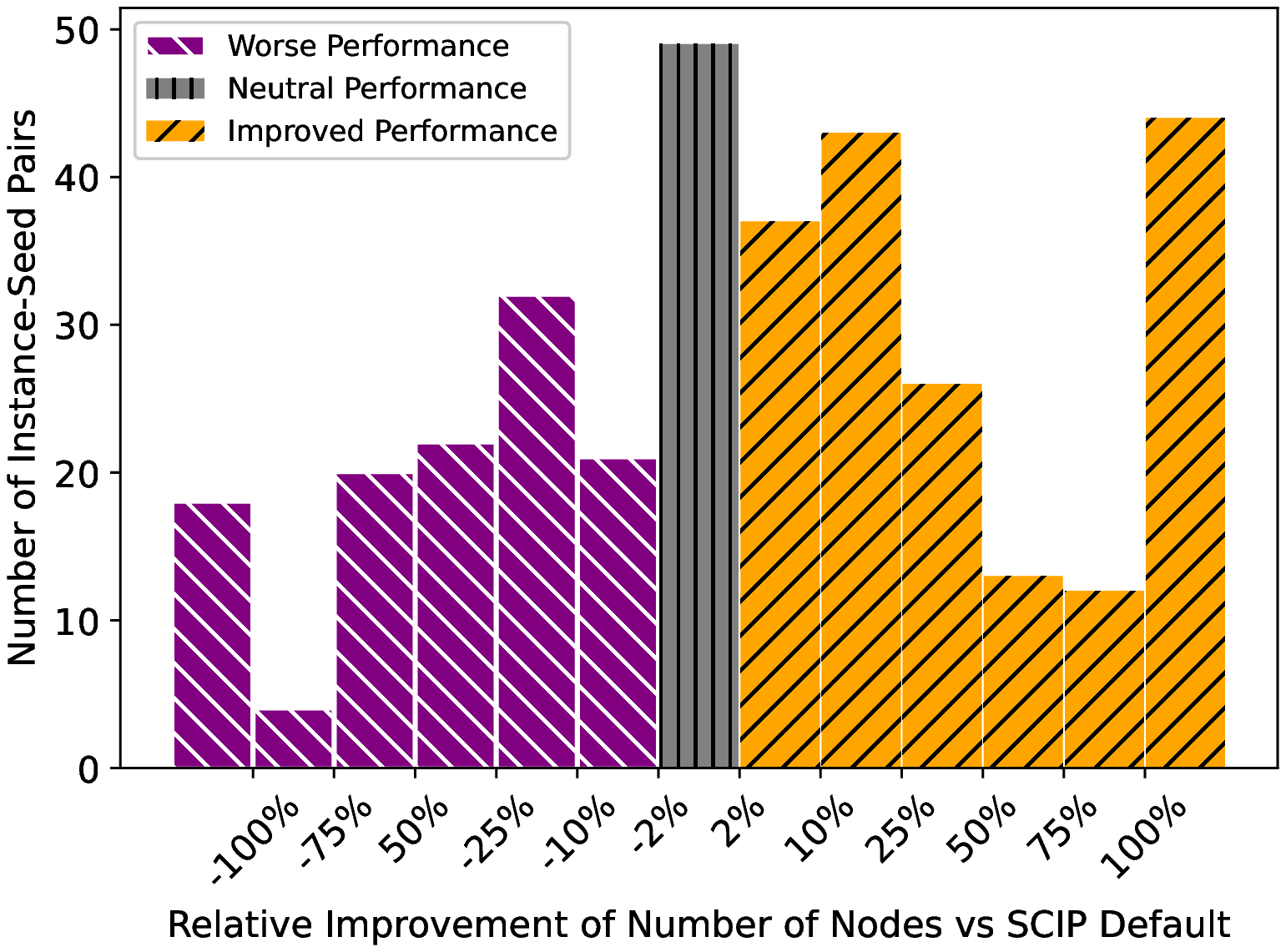}\hspace{.1em}%
    \includegraphics[height=6.5cm,width=0.48\linewidth]{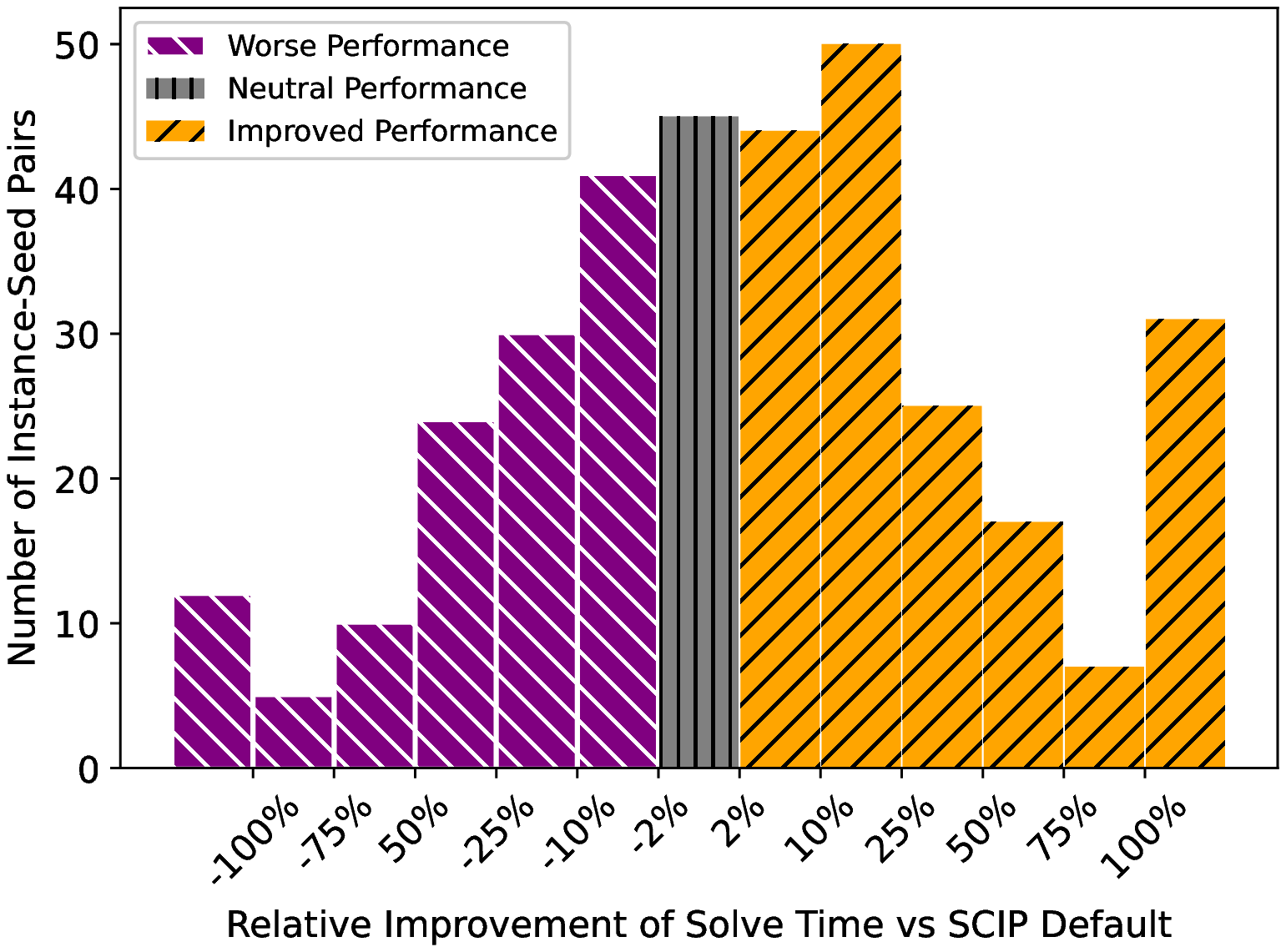}
    \caption{MIPLIB with an optimal solution provided and heuristics disabled.}
    \end{subfigure}
    \caption{Bar plots of relative improvement of \emph{gmi-}$10^{-5}$ compared to default SCIP over affected instance-seed pairs. (Left) Number of nodes. (Right) Solve time.}
    \label{fig:bar_plot_miplib}
\end{figure}

Figure~\ref{fig:bar_plot_miplib} shows the distribution of performance improvement per affected instance-seed pair. We observe a diverse distribution, with the majority of instance-seed pairs either performing at least $10\%$ better or at least $10\%$ worse in both number of nodes and solve time. This is surprising as our method only adds a relatively small value to the weighted sum branching score, and therefore acts as a tiebreaker. While many instances exhibit worse performance on \emph{gmi-}$10^{-5}$, there are consistently more instance-seed pairs that perform correspondingly better than default SCIP over all levels of improvement. This is evident both for number of nodes and for solve time, and is confirmed by individual Wilcoxon signed-rank tests over all data sets, where a p-value of at most 0.028 was observed. 

\section{Conclusion}

In this paper, we developed a new branching rule based on the correspondence between Gomory mixed-integer cuts and split disjunctions, leveraging the efficacy of cutting planes as a measure for the relevance of a variable for branching. We used the branching rule with both unstrengthened and strengthened cuts, showing that the unstrengthened versions, which are directly derived from potential branching decisions, provide a better measure to reduce the number of nodes, and solve time. Our branching rule results in low numbers of nodes while being less expensive than strong branching. When integrated in the state-of-the-art hybrid branching algorithm of SCIP, the score provided by our branching rule significantly reduces both solve time and number of nodes over MIPLIB 2017. Future work includes extending our idea beyond Gomory mixed-integer cuts, to any cut that is linked to a split disjunction, e.g., lift-and-project cuts.

\section*{Acknowledgements}
We thank Tobias Achtherberg, Antonia Chmiela, and Leona Gottwald for insightful discussions that helped this paper. The work for this article has been conducted in the Research Campus MODAL funded by the German Federal Ministry of Education and Research (BMBF) (fund numbers 05M14ZAM, 05M20ZBM). The described research activities are funded by the Federal Ministry for Economic Affairs and Energy within the project UNSEEN (ID: 03EI1004-C).

\bibliographystyle{unsrt}
\bibliography{mybib}

\end{document}